\documentclass[11pt]{amsart}
\usepackage{amscd,graphicx,psfrag}

\newcommand{\p}{\partial}
\newcommand{\Z}{\mathbb Z}
\newcommand{\Q}{\mathbb Q}
\newcommand{\G}{\mathcal G}

\newcommand{\D}{\mathcal D}
\newcommand{\B}{\mathcal B}
\newcommand{\C}{\mathcal C}
\newcommand{\M}{\mathcal M}
\newcommand{\R}{\mathcal R}
\renewcommand{\S}{\mathcal S}

\newcommand{\rhofo}{\rho}

\newcommand{\im}{\operatorname{im}}

\renewcommand{\phi}{\varphi}

\newcommand{\tr}{\operatorname{tr}}

\newcommand{\su}{\mathfrak{su}}

\newcommand{\ad}{\operatorname{ad}}
\newcommand{\rk}{\operatorname{rk}}

\newcommand{\Hess}{\operatorname{Hess}}

\newcommand{\Lef}{\operatorname{Lef}\,}
\newcommand{\hol}{\operatorname{hol}}
\newcommand{\arf}{\operatorname{arf}}
\newcommand{\sign}{\operatorname{sign}}
\newcommand{\Span}{\operatorname{Span}}
\newcommand{\ZZ}{\Z[\Z]} % not {\Z_2\oplus \Z_2} !

\renewcommand{\P}{\mathcal{PR}}
\newcommand{\Id}{\operatorname{Id}}%

\newtheorem{theorem}{Theorem}[section]
\newtheorem{lemma}[theorem]{Lemma}
\newtheorem{conjecture}[theorem]{Conjecture}
\newtheorem{proposition}[theorem]{Proposition}
\newtheorem{corollary}[theorem]{Corollary}%

\def\acknowledgmentname{Acknowledgments.}
\newenvironment{acknowledgment}
   {
    \vspace*{\baselineskip}
    \noindent{\small\textbf{\acknowledgmentname}}
    \unskip\noindent}{}

\begin{document}%

\title{Casson--type invariants in dimension four}%
\thanks{The first author was partially supported by NSF Grants
9971802 and 0204386. The second author was partially supported by NSF
Grant 0305946 and Max Planck Institut f\"ur Mathematik in Bonn}
\author[Daniel Ruberman]{Daniel Ruberman}
\address{Department of Mathematics,
MS 050,
Brandeis University,
Waltham, MA 02454}
\email{\rm{ruberman@brandeis.edu}}
\author[Nikolai Saveliev]{Nikolai Saveliev}
\address{Department of Mathematics,
University of Miami,
Box 249085,
Coral Gables, FL 33124}
\email{\rm{saveliev@math.miami.edu}}

%\keywords{Casson invariant, Rohlin invariant, Floer homology, flat
%moduli spaces}


%\begin{abstract}
%\end{abstract}

\maketitle%

\section{Introduction}
This article surveys our ongoing project about the relationship
between invariants extending the classical Rohlin invariant of
homology spheres and those coming from 4--dimensional (Yang-Mills)
gauge theory; it will appear in the Proceedings of the Fields-McMaster Conference on Geometry and Topology of Manifolds.  We are mainly concerned with a special class of
manifolds for which the two types of invariants are defined and
can be compared. This class contains, in particular, manifolds
having the homology of $S^1\times S^3$.  Rohlin's theorem about
the signature of closed smooth spin 4--manifolds gives rise to a
mod--2 invariant of a homology $S^1\times S^3$. On the gauge
theoretic side, the invariant is obtained by counting flat
connections on appropriate bundles. This count is inspired by
Donaldson's~\cite{donaldson-kronheimer} count of anti-self-dual
connections on $SU(2)$ bundles; its flat analogue was first
studied by Furuta and Ohta~\cite{furuta-ohta}.

The main conjecture towards which this project is directed is that
the Rohlin invariant and the gauge theoretic invariant coincide 
for homology $S^1\times S^3$. The model for the whole discussion 
is Casson's beautiful theorem relating his invariant (in its gauge 
theoretic manifestation as described by Taubes~\cite{taubes:casson}) 
and Rohlin's invariant of homology 3--spheres. We will discuss the 
implications of this conjecture for some classical problems in 
low-dimensional topology, and progress we have made towards proving 
the conjecture.  This progress includes the verification of the 
conjecture in some special cases, a `surgery' program for approaching 
the conjecture by expanding its scope to include a wider category of 
manifolds, and the verification of this expanded conjecture for 
homology 4--tori. Much of this material is contained in our three 
papers \cite{ruberman-saveliev:casson,ruberman-saveliev:mappingtori,
ruberman-saveliev:donaldson} but we have included a broader overview 
as well as some additional examples.

\begin{acknowledgment} We would like to thank Scott Baldridge for pointing 
out the manifolds described in Section~\ref{S:circle}, and Liviu Nicolaescu 
for his input on computing orientations of flat moduli spaces over these 
manifolds. We also thank Kim Fr\o yshov for some interesting conversations
related to the material presented in this paper.  We express our appreciation 
to the conference organizers  for providing such a stimulating environment at 
McMaster.
\end{acknowledgment}

%%%%%%%%%%%%%%%%%%%%%%%%%%%%%%%%%%%%%%%%%%%%%%%%%%%%%%%%%%%%%%%%%%%%%%%%%

\section{The Rohlin invariant}\label{S:rohlin}
This section is a review of the homology cobordism group and the Rohlin
invariant and some of their applications in topology.  More
information and references may be found in the
books~\cite{saveliev:casson,saveliev:spheres}.

%%%%%%%%%%%%%%%%%%%%%%%%%%%%%%%%%%%%%%%%%%%%%%%%%%%%%%%%%%%%%%%%%%%%%%%%%

\subsection{Homology spheres}
By an (integral) homology sphere we will mean a closed oriented 3--manifold
$\Sigma$ such that $H_*(\Sigma;\Z) = H_*(S^3;\Z)$.  According to the
Poincar\'e conjecture, every simply connected homology sphere is homeomorphic
to $S^3$. There exist, however, many non-simply connected homology spheres,
as the following examples demonstrate.

For any three positive pairwise relatively prime integers $p$, $q$, and $r$,
the zero set of the complex polynomial $x^p + y^q + z^r$ is a complex surface
which has an isolated singularity at the origin. The link of this singularity,
\[
\Sigma(p,q,r) = \{\,x^p + y^q + z^r = 0\,\}\,\cap\,S^5,
\]
is a homology sphere sphere referred to as a Brieskorn homology sphere.  The
Brieskorn homology sphere $\Sigma(2,3,5)$ is also known as the Poincar\'e 
sphere. The above construction can be generalized to obtain Seifert fibered 
homology spheres $\Sigma(a_1,\ldots,a_n)$  for any positive pairwise 
relatively prime integers $a_1,\ldots,a_n$.

A more general construction of homology spheres is as follows.  Let $k$ be a
knot in $S^3$ and $q$ an integer then any manifold $S^3 + (1/q)\,k$ obtained
by $(1/q)$-surgery of $\Sigma$ along $k$ is a homology sphere.  For example,
$\Sigma(2,3,5)$ can be obtained by $(-1)$--surgery on the left handed trefoil.
Not all homology spheres can be obtained by this construction;  however, if
one allows surgery along links of more than one component,  one obtains all
homology spheres (and in fact all closed oriented 3--manifolds).

%%%%%%%%%%%%%%%%%%%%%%%%%%%%%%%%%%%%%%%%%%%%%%%%%%%%%%%%%%%%%%%%%%%%%%%%%%%

\subsection{Homology cobordisms}
A homology cobordism from a homology sphere $\Sigma_0$ to a homology sphere 
$\Sigma_1$ is a smooth compact oriented 4--manifold $W$ with boundary $\p W 
= - \Sigma_0 \cup \Sigma_1$ such that the inclusions $\Sigma_i\to W$ induce 
isomorphisms $H_*(\Sigma_i;\Z)\to H_*(W;\Z)$ for $i = 0,1$.

An obvious example of a homology cobordism from a homology sphere $\Sigma$
to itself is the product $W = [0,1] \times \Sigma$. The mapping cylinder
$W_{\tau}$ of any orientation preserving diffeomorphism $\tau: \Sigma \to
\Sigma$ is also a homology cobordism from $\Sigma$ to itself;  note that
$W_{\tau}$ is diffeomorphic to the product $[0,1]\times \Sigma$  but need
not be diffeomorphic to it rel boundary.

The first example of a non-trivial homology cobordism was constructed by
Mazur; his manifold is a simply connected homology cobordism between $S^3$
and $\Sigma(2,5,7)$. Many more examples of homology cobordisms can be
obtained by performing surgery along knot and link concordances in $[0,1]
\times S^3$, and by various other constructions.

%%%%%%%%%%%%%%%%%%%%%%%%%%%%%%%%%%%%%%%%%%%%%%%%%%%%%%%%%%%%%%%%%%%%%%%%%%%%

\subsection{Homology cobordism group}
The relation of being homology cobordant is an equivalence relation on the
class of homology spheres. The set of equivalence classes of homology
spheres with the operation induced by connected sum is an abelian group
called the homology cobordism group and denoted $\Theta^3$. The inverse
element of $[\Sigma]\in \Theta^3$ is $[-\Sigma]$, and the zero element is
the homology cobordism class of $S^3$.

%%%%%%%%%%%%%%%%%%%%%%%%%%%%%%%%%%%%%%%%%%%%%%%%%%%%%%%%%%%%%%%%%%%%%%%%%%%%

\subsection{The Rohlin invariant}
Every homology sphere $\Sigma$ is the boundary of a smooth compact spin
4--manifold $X$, whose signature is necessarily divisible by eight. Additivity
of the signature and the Rohlin theorem  (which asserts that the signatures
of any two such manifolds bounding $\Sigma$ can only differ by a multiple of 
sixteen) imply that the quantity
\[
\rho(\Sigma) = \frac 1 8\,\sign X \pmod 2
\]
only depends on $\Sigma$ and takes values in $\Z_2 = \Z/2\Z$. It is called
the Rohlin invariant of $\Sigma$. One can easily see that $\rho$ defines a
homomorphism $\rho: \Theta^3 \to \Z_2$.

As an example, consider the singularity at zero of $x^2 + y^3 + z^5 = 0$. It
has a resolution $X$ which is a smooth compact simply connected spin manifold
with boundary $\Sigma(2,3,5)$ and intersection form $E_8$. Therefore,
$\sign X = -8$ and $\rho (\Sigma(2,3,5)) = 1\pmod 2$. In particular, we see
that $\rho: \Theta^3 \to \Z_2$ is an epimorphism.

The Rohlin invariant can be defined more generally for pairs $(Y,\sigma)$,
where $Y$ is a closed oriented 3-manifold and $\sigma$ is a spin structure
on $Y$. By definition,
\[
\rho(Y,\sigma) = \frac 1 8\,\sign X\pmod 2,
\]
where $X$ is any smooth compact spin 4--manifold with (spin) boundary
$(Y,\sigma)$. This invariant takes values in $\Q/2\Z$. A homology sphere
$\Sigma$ has a unique spin structure hence it does not show up in the
notation $\rho(\Sigma)$.

%%%%%%%%%%%%%%%%%%%%%%%%%%%%%%%%%%%%%%%%%%%%%%%%%%%%%%%%%%%%%%%%%%%%%%%%%%%

\subsection{The structure of $\Theta^3$}
The fact that $\rho: \Theta^3 \to \Z_2$ is an epimorphism had been the only
known general fact about $\rho$ and $\Theta^3$ until the early 1980's when
some progress was made using gauge theory. The Donaldson diagonalizability
theorem for the intersection forms of smooth closed oriented definite
4--manifolds implies right away that $\Theta^3$ is an infinite group.
In fact, one can see that $\Sigma(2,3,5)$ is an element of infinite order
in $\Theta^3$. More elements of infinite order were found by Fintushel
and Stern using equivariant gauge theory. Later Furuta showed that
$\Theta^3$ is infinitely generated -- in fact,  Brieskorn homology spheres
$\Sigma(2,3,6m-1)$ with $m\ge 1$ all have infinite order in $\Theta^3$ and
are linearly independent over $\Z$.

The question of whether $\Theta^3$ has torsion remains open; the latest
result in this direction is as follows, see \cite{fukumoto-furuta},
\cite{fukumoto-furuta-ue} and \cite{saveliev:fukumoto-furuta}.

\begin{theorem}\label{T:ffs}
Let $\Sigma$ be a homology sphere which is homology cobordant to a Seifert
fibered homology sphere. If $\rho(\Sigma)$ is non-trivial then $\Sigma$ has
infinite order in $\Theta^3$.
\end{theorem}

It is not known if all homology spheres are homology cobordant to
Seifert fibered ones but this is considered highly unlikely.  In
particular, Fr\o yshov has an (unpublished) extension of his work on the
$h$--invariant \cite{froyshov}, which seems to give rise to
counterexamples.  

Describing the structure of $\Theta^3$ is interesting in its own right but
also because of its applications some of which are described below, see
also Section \ref{S:surgery-4}.

%%%%%%%%%%%%%%%%%%%%%%%%%%%%%%%%%%%%%%%%%%%%%%%%%%%%%%%%%%%%%%%%%%%%%%%%

\subsection{Triangulation conjecture}\label{S:rohlin-tri}
The triangulation conjecture in dimension $n$ asserts that every (closed)
topological $n$--manifold is homeomorphic to a simplicial complex.  This
conjecture has been long known to hold in dimensions three and lower. It
fails in dimension four, which follows by combining Freedman's
classification of simply connected topological manifolds and the Casson
invariant theory, see Section \ref{S:casson-tri}.

The triangulation conjecture remains open in dimensions $n\ge 5$. Amazingly
enough, it is equivalent in these dimensions to a specific question about
the structure of $\Theta^3$. Namely, according to a theorem of
Matumoto~\cite{matumoto:triangulation} and
Galewski and Stern~\cite{galewski-stern:simplicial}, the
triangulation conjecture holds in all dimensions
$n \ge 5$ if and only if there exists a homology sphere $\Sigma$ of order
two in $\Theta^3$ having non-trivial Rohlin invariant.  Needless to say, 
no such homology sphere has been found; Theorem \ref{T:ffs} implies that 
there is no need to search among Seifert fibered homology spheres and the 
homology spheres homology cobordant to them.

%%%%%%%%%%%%%%%%%%%%%%%%%%%%%%%%%%%%%%%%%%%%%%%%%%%%%%%%%%%%%%%%%%%%%%%%%

\subsection{Simply connected homology cobordisms}\label{S:rohlin-sc}
Every homology sphere $\Sigma$ is homology cobordant to itself via the
product homology cobordism $[0,1]\times \Sigma$. This cobordism is not
simply connected unless $\Sigma$ is. A natural question arises whether
$\Sigma$ can be homology cobordant to itself via a \emph{simply connected}
homology cobordism. For some homology spheres the answer to this question
is positive, while for others it is negative. For example, there are no
simply connected homology cobordisms of the Poincar\'e homology sphere to
\emph{any} homology sphere, not just to itself, see Taubes
\cite{taubes:periodic}. More examples of such behavior have been given by
Fintushel and Stern \cite{fs:instanton}.

It has been conjectured that no two homology spheres with non-trivial
Rohlin invariant can be homology cobordant via a simply
connected homology cobordism (note that the examples of Fintushel and
Stern may have both trivial and non-trivial Rohlin invariants).

%%%%%%%%%%%%%%%%%%%%%%%%%%%%%%%%%%%%%%%%%%%%%%%%%%%%%%%%%%%%%%%%%%%%%%%%%

\section{The Casson invariant}\label{S:casson}
The Casson invariant is an integer valued invariant of homology spheres,
defined by Casson in 1985. For a homology sphere $\Sigma$, it was defined
as a `creative' count of $SU(2)$ representations of $\pi_1\Sigma$, see
Akbulut-McCarthy \cite{akbulut-mccarthy}. Later, Taubes
\cite{taubes:casson} reformulated this original definition in
gauge theoretic terms using the fact that every $SU(2)$ representation of
$\pi_1\Sigma$ arises as the holonomy representation of a flat connection.

%%%%%%%%%%%%%%%%%%%%%%%%%%%%%%%%%%%%%%%%%%%%%%%%%%%%%%%%%%%%%%%%%%%%%%%%%

\subsection{Definition of the Casson invariant}\label{S:casson-def}
Let $P \to \Sigma$ be a trivialized $SU(2)$ bundle over a homology sphere
$\Sigma$ and consider the moduli space $\R^*(\Sigma)$ of (gauge equivalence
classes) of irreducible flat connections on $P$. After perturbation of the
flatness equation $F_A = 0$, if necessary, the moduli space $\R^*(\Sigma)$
is a compact oriented zero-dimensional manifold. Define the Casson invariant
\begin{equation}\label{E:casson}
\lambda(\Sigma) = \frac 1 2\,\#\R^*(\Sigma),
\end{equation}
where $\#\R^*(\Sigma)$ stands for the signed count of the (finitely many)
points in $\R^*(\Sigma)$. Making sense of this definition requires a lot
of work, which goes mostly into orienting $\R^*(\Sigma)$ and choosing proper
perturbations.

%%%%%%%%%%%%%%%%%%%%%%%%%%%%%%%%%%%%%%%%%%%%%%%%%%%%%%%%%%%%%%%%%%%%%%%%%%%

\subsection{Orientation}\label{S:casson-or}
We say that $A \in \R^*(\Sigma)$ is non-degenerate if $H^1(\Sigma;\ad A) =
0$. This is equivalent to saying that, for a choice of Riemannian metric
on $\Sigma$, the elliptic operator
\begin{equation}\label{E:ka}
K_A = \begin{pmatrix} 0 & d_A^* \\ d_A & - *d_A
\end{pmatrix}
\end{equation}
acting on the space $(\Omega^0\,\oplus\,\Omega^1)(\Sigma;\ad P)$ has zero
kernel. If all points in $\R^*(\Sigma)$ are non-degenerate then there are
only finitely many of them. The sign of $A\in \R^*(\Sigma)$ is defined as
$(-1)^{\mu(A)}$ where $\mu(A)$ is the spectral flow of the family of
operators $K_{A(t)}$ along a path $A(t)$ from the product connection
$\theta$ to $A$. The quantity $\mu(A)$ reduced modulo 8 only depends on
$A$ and not on a particular choice of $A(t)$, and is referred to as the
Floer index of $A$. If $\R^*(\Sigma)$ fails to be non-degenerate, it is
perturbed and then oriented using a similar procedure.

%%%%%%%%%%%%%%%%%%%%%%%%%%%%%%%%%%%%%%%%%%%%%%%%%%%%%%%%%%%%%%%%%%%%%%%%%%%

\subsection{Perturbations}\label{S:casson-pert}
In the degenerate situation,  the flatness equation $F_A = 0$ is perturbed
into $F_A = *\nabla h$, where $h$ is a function on the connections defined
as follows, see Taubes \cite{taubes:casson} and Herald
\cite{herald:perturbations}.

A collection of disjoint closed loops $\gamma_k$ embedded in $\Sigma$ will
be called a link. Given a link of $n$ loops, consider $n$ smooth functions
$f_k: SU(2)\to\mathbb R$ invariant with respect to conjugation, and define
$h(A)$ as the sum of $f_k(\hol_A(\gamma_k))$. For analytical reasons, one
uses a little more sophisticated definition of $h(A)$ obtained by averaging
the above over the neighboring links. More precisely, thicken the link
$\{\,\gamma_k\,\}$ into a collection of embeddings $\gamma_k: S^1\times D^2 
\to \Sigma$ with disjoint images, and define
\[
h(A) = \sum_{k=1}^n\;\int_{D^2}\;f_k(\hol_A(\gamma_k(S^1\times\{z\})))\,
\eta(z)\,d^2 z,
\]
where $\hol_A(\gamma_k(S^1\times\{z\}))$ stands for holonomy of $A$ around
the loop $\gamma_k(S^1 \times\{z\})$, $z\in D^2$, and $\eta$ is any smooth
rotationally symmetric bump function on $D^2$ with support away from the
boundary of $D^2$ and with integral one.

The moduli space of irreducible solutions of the equation $F_A = *\nabla h$
is denoted by $\R^*_h (\Sigma)$; a point $A$ in this space is called
non-degenerate if the kernel of the operator $K_A + \Hess_A h$ vanishes. It
turns out that there exist the so called abundant links which have the
property that, for any small generic functions $f_k$, the moduli space
$\R^*_h(\Sigma)$ is non-degenerate. The non-degenerate space $\R^*_h(\Sigma)$
consists of finitely many points. We orient it using the spectral flow of the
operators $K_{A(t)} + \Hess_{A(t)} h$ as above.

The class of perturbations described above is large enough for our current
purposes. We expand this class to include loops with a common base point to
deal with equivariant aspects of gauge theory in Sections \ref{S:free} and
\ref{S:surgery-3}, compare with Floer \cite{floer:instanton} and Herald 
\cite{herald:basepoint}.

%%%%%%%%%%%%%%%%%%%%%%%%%%%%%%%%%%%%%%%%%%%%%%%%%%%%%%%%%%%%%%%%%%%%%%%%%%%%%

\subsection{Properties of the Casson invariant}\label{S:casson-prop}
The above construction ends up in an integer valued invariant  $\lambda
(\Sigma)$ which only depends on $\Sigma$ and not on the choices made in
its definition. We will concentrate on one specific property of this
invariant, namely, that $\lambda(\Sigma) = \rho(\Sigma)\pmod 2$. To verify
this property, one traditionally goes back to the original definition of
the Casson invariant using representations spaces. One proves a surgery
formula for the Casson invariant in terms of the Alexander polynomial,
and then shows that the modulo 2 reduction of this formula gives the
surgery formula for the Rohlin invariant. A more direct proof using the
gauge theoretic definition of the Casson invariant and avoiding the
Alexander polynomial altogether will be given in Section \ref{S:surgery-3}.

%%%%%%%%%%%%%%%%%%%%%%%%%%%%%%%%%%%%%%%%%%%%%%%%%%%%%%%%%%%%%%%%%%%%%%%%%%%%%

\subsection{Triangulation of manifolds in dimension four}\label{S:casson-tri}
The existence of an invariant $\lambda(\Sigma)$ lifting the Rohlin invariant
to the integers provides for a negative solution of the triangulation
conjecture in dimension four. Here is a sketch of the argument. 

Let $X$ be a closed topological 4--manifold with intersection form $E_8$
(such a manifold exists by Freedman's classification), and suppose that it
is homeomorphic to a simplicial complex. Without loss of generality, we
will assume that the links of all vertices in this complex but maybe one
are homeomorphic to $S^3$. The link of the remaining vertex, say $v$, is
then a homotopy sphere $\Sigma$. Since $\pi_1\Sigma$ is trivial, it has no
irreducible $SU(2)$ representations and hence $\rho(\Sigma) = \lambda
(\Sigma) = 0\pmod 2$. On the other hand, removing an open neighborhood of
$v$ makes $X$ into a piecewise linear manifold (which is the same as a
smooth manifold in dimension four) with boundary $\Sigma$ and the
intersection form $E_8$. Hence $\rho(\Sigma) = 1\pmod 2$, and this
contradiction shows that $X$ could not have been homeomorphic to a
simplicial complex.

The only use for the Casson invariant in this application is that it shows
that the Rohlin invariant of a homotopy sphere vanishes. This fact can also
be deduced from the Poincar\'e conjecture.

%%%%%%%%%%%%%%%%%%%%%%%%%%%%%%%%%%%%%%%%%%%%%%%%%%%%%%%%%%%%%%%%%%%%%%%%%%%%%

\section{The Furuta--Ohta invariant}\label{S:fo}
Both Casson and Rohlin invariant are generalized in this section to invariants
of $\ZZ$--homology $S^1\times S^3$, which are smooth 4-manifolds satisfying
certain homological conditions. The equality $\lambda(\Sigma) = \rho(\Sigma)
\pmod 2$ that we had in dimension three becomes a conjecture for these
generalized invariants. We explain how this conjecture is related to the
triangulation conjecture in dimensions five and higher.

%%%%%%%%%%%%%%%%%%%%%%%%%%%%%%%%%%%%%%%%%%%%%%%%%%%%%%%%%%%%%%%%%%%%%%%%%%%%%

\subsection{Homology $S^1 \times S^3$}\label{S:fo-ex}
A $\ZZ$--homology $S^1 \times S^3$ is a smooth oriented 4-manifold $X$ such
that $H_*(X;\Z) = H_* (S^1\times S^3;\Z)$ and $H_*(\tilde X;\Z) =
H_*(S^3;\Z)$, where $\tilde X$ is the universal abelian cover of $X$.

An ample source of examples of $\ZZ$--homology $S^1\times S^3$ is provided 
by the following
operation. Given a homology cobordism $W$ from a homology sphere $\Sigma$
to itself, form its closure $\overline W$ by identifying the two copies
of $\Sigma$ in the boundary of $W$ by the identity map. The closure
$\overline W$ is always a $\ZZ$--homology $S^1\times S^3$. For instance, the
closure of the product cobordism $[0,1] \times \Sigma$ is the product
$S^1\times \Sigma$, and the closure of a mapping cylinder $W_{\tau}$ is the
mapping torus $X_{\tau} = ([0,1] \times \Sigma)/(0,x)\sim (1,\tau(x))$.

Another series of examples is generated by $S^1$--bundles $X \to Y$ over
3--manifolds $Y$ with $H_*(Y;\Z) = H_*(S^1\times S^2;\Z)$. These
examples were pointed out to us at the McMaster conference by Scott
Baldridge, as they arise in his work on circle actions and
Seiberg-Witten theory.  An application
of the Gysin exact sequence shows that, if the Euler class $e\in H^2(Y;\Z)$
of the $S^1$--bundle generates the group $H^2 (Y;\Z) = \Z$, the manifold
$X$ has integral homology of $S^1\times S^3$. However, the second condition
that $H_*(\tilde X;\Z) = H_*(S^3;\Z)$ is only satisfied if $H_*(\tilde Y;\Z)
= H_* (S^2;\Z)$, which is equivalent to saying that the Alexander polynomial
of $Y$ is trivial.  It should be noted that there are plenty of manifolds
$Y$ with homology of $S^1\times S^2$ whose Alexander polynomial is trivial;
for instance, any $Y$ obtained by 0--surgery on a knot in $S^3$ with trivial
Alexander polynomial (such as untwisted Whitehead double) will do.

More examples come from 2--knots in the 4--sphere, surgery along which 
produces a manifold with the homology  of $S^1 \times S^3$. To satisfy the 
homological condition on the universal abelian cover of this manifold, one
would start with a 2--knot with trivial Alexander polynomial. Such knots are
readily constructed, for example as the $k$--fold twist-spin of a knot whose
$k$--fold branched cover is a homology sphere.

%%%%%%%%%%%%%%%%%%%%%%%%%%%%%%%%%%%%%%%%%%%%%%%%%%%%%%%%%%%%%%%%%%%%%%%%%%%%%

\subsection{Extension of the Rohlin invariant}\label{S:fo-rohlin}
Let $X$ be a $\ZZ$--homology $S^1\times S^3$ and choose an embedded
3--manifold $M\subset X$ whose fundamental class generates $H_3(X;\Z) = \Z$.
Note that $M$ need not be a homology sphere. We define
\[
\rhofo(X) = \rho(M,\sigma)\pmod 2
\]
where $\sigma$ is a spin structure on $M$ induced from $X$. This is a well
defined invariant of $X$ independent of the choices made in its definition,
compare with \cite{ruberman:ds}. For example, if $X = \overline W$ is the
closure of a homology cobordism $W$ from a homology sphere $\Sigma$ to
itself then $\rhofo(X) = \rho(\Sigma)\pmod 2$.

%%%%%%%%%%%%%%%%%%%%%%%%%%%%%%%%%%%%%%%%%%%%%%%%%%%%%%%%%%%%%%%%%%%%%%%%%%%%%

\subsection{Definition of the Furuta--Ohta invariant}\label{S:fo-def}
Let $X$ be a $\ZZ$--homology $S^1\times S^3$ and $\M^*(X)$ the moduli space
of irreducible ASD connections on a trivial $SU(2)$ bundle $P \to X$.  Note
that all such connections are flat by Chern--Weil theory.

The formal dimension of $\M^*(X)$ is equal to $8\,c_2(P)\,[X] - 3(1 - b_1 +
b^+_2)(X) = 0$. After perturbing the ASD equation $F^+_A = 0$ if necessary,
$\M^*(X)$ is a compact oriented zero dimensional manifold (it
is worth mentioning that the perturbed ASD connections need no longer be
flat or even perturbed flat). The Furuta--Ohta invariant of $X$ is then
defined as
\begin{equation}\label{E:fo}
\lambda_{FO} (X) = \frac 1 4\,\#\M^*(X),
\end{equation}
where $\#\,\M^*(X)$ stands for a signed count of points in $\M^*(X)$. The
compactness of $\M^*(X)$ is guaranteed by the condition $H_*(\tilde X;\Z)
= H_*(S^3;\Z)$, see Furuta--Ohta \cite{furuta-ohta}. As with the Casson
invariant, the main work goes into orienting $\M^*(X)$ and choosing proper
perturbations.

%%%%%%%%%%%%%%%%%%%%%%%%%%%%%%%%%%%%%%%%%%%%%%%%%%%%%%%%%%%%%%%%%%%%%%%%%%%%%

\subsection{Orientation}\label{S:fo-or}
The moduli space $\M^*(X)$ is called non-degenerate if the ASD operator
\begin{equation}\label{E:asd}
D_A = d_A^*\,\oplus\,d^+_A: \Omega^1(X;\ad P) \to (\Omega^0\,\oplus\,
\Omega^2_+) (X;\ad P)
\end{equation}
has trivial cokernel for every $A \in \M^*(X)$. A non-degenerate $\M^*(X)$
is oriented using the following construction from the Donaldson theory. Let
$\B(X)$ be the space of the gauge equivalence classes of connections on $P
\to X$, and $\Lambda_X$ the determinant bundle of the family $D_A$ over
$\B(X)$. This is a real line bundle with the property that, over $\M^*(X)
\subset \B(X)$, it restricts to the orientation bundle of $\M^*(X)$.
According to Donaldson, the bundle $\Lambda_X$ is trivial over $\B(X)$, and
a choice of trivialization of $\Lambda_X$ given by an orientation of
$H^1(X;\mathbb R) = \mathbb R$ (called homology orientation) fixes an
orientation of $\M^*(X)$. If $\M^*(X)$ fails to be non-degenerate, it is
perturbed first and then oriented using a similar construction.

%%%%%%%%%%%%%%%%%%%%%%%%%%%%%%%%%%%%%%%%%%%%%%%%%%%%%%%%%%%%%%%%%%%%%%%%%%%%

\subsection{Perturbations}\label{S:fo-pert}
To make sense of the above definition of $\lambda_{FO}(X)$,  the ASD equation
$F^+_A = 0$ defining $\M^*(X)$ may need to be perturbed into $F^+_A = 
\sigma(A)$ using admissible perturbations $\sigma(A) \in \Omega^2_+(X;\ad P)$. 
The latter are constructed as follows, compare with Donaldson 
\cite{donaldson:orientation}.

Let us consider an embedding $\psi: S^1 \to X$ and extend it to an embedding
$\psi: S^1\times N^3 \to X$ where $N^3$ is an oriented 3--manifold. For any
connection $A$ in $P$ denote by $\hol_A (\psi (S^1\times\{x\}),s)\in SU(2)$
the holonomy of $A$ around the loop $\psi (S^1\times\{x\})$ starting at the
point $\psi (s,x)$. Let $\Pi: SU(2)\to \su(2)$ be the projection given by
\[
\Pi (u) = u - \frac 1 2\,\tr(u)\cdot\Id.
\]
Assigning $\Pi\hol_A(\psi(S^1\times\{x\}),s)$ to $\psi(s,x)\in X$ defines
a section of $\ad P$ over $\psi (S^1\times N^3)$.  Now, given a form $\nu
\in \Omega^2_+(X)$ supported in $\psi(S^1\times N^3)$, define a section
\[
\sigma (\nu, \psi, A) \in \Omega^2_+ (X,\ad P)
\]
by taking tensor product of $\Pi \hol_A (\psi (S^1\times\{x\}),s)$ with
$\nu$ over $\psi (S^1\times N^3)$ and letting it be zero otherwise.

More generally, consider a collection of embeddings $\psi_k: S^1\to X$,
$k = 1,\ldots, n$, with disjoint images, called a link, and extend it to a
collection of embeddings $\psi_k: S^1\times N^3_k\to X$ as above so that
the $\psi_k (S^1\times N^3_k)$ are still disjoint.  For any choice of $n$
smooth functions $\bar f_1,\ldots,\bar f_n: [-2,2]\to\mathbb R$ with
vanishing derivatives at $\pm\,2$, define admissible perturbation
\[
\sigma (A) = \sum_{k=1}^n\; \p \bar f_k\cdot\sigma(\nu_k, \psi_k, A),
\]
where $\p \bar f_k$ is the function $\bar f'_k$ evaluated at $\tr\hol_A
(\psi_k (S^1\times\{x\}),s)$, and $\nu_k$ are real valued self--dual forms
on $X$, each supported in its respective $\psi_k(S^1\times N^3_k)$.

Given an admissible perturbation $\sigma$, the set of the gauge equivalence
classes of irreducible solutions of the equation $F^+_A = \sigma(A)$ will
be denoted by $\M^*_{\sigma}(X)$. There exist the so called abundant links
which have the property that, for any small generic functions $\bar f_k$,
the moduli space $\M^*_{\sigma}(X)$ is non-degenerate (the non-degeneracy
condition here means that the cokernel of the perturbed ASD operator
(\ref{E:asd}) vanishes). Once $\M^*_{\sigma}(X)$ is non-degenerate it is 
oriented as in Section \ref{S:fo-or}.

The class of admissible perturbations will be expanded in Section
\ref{S:tori-pert} to handle the equivariant gauge theory.

%%%%%%%%%%%%%%%%%%%%%%%%%%%%%%%%%%%%%%%%%%%%%%%%%%%%%%%%%%%%%%%%%%%%%%%%%%%

\subsection{Properties of the Furuta--Ohta invariant}\label{S:fo-prop}
One can show that the above construction ends up in an invariant which only
depends on $X$ and a choice of orientation of $H^1(X;\mathbb R)$. Note that
this invariant can be viewed as one quarter of a degree zero Donaldson
polynomial of $X$, except the latter is formally not defined for trivial
bundles or manifolds with $b^+_2 = 0$.

Next we wish to discuss the factor of one quarter in the definition of the
Furuta--Ohta invariant (\ref{E:fo}) vs. the factor of one half for the 
Casson invariant (\ref{E:casson}). The reason for the extra one half is a 
2--fold symmetry arising from the action of $H^1(X;\Z_2) = \Z_2$ on 
$\M^*(X)$. This action can be described as follows.

Let us view $\chi \in H^1(X;\Z_2)$ as a homomorphism from $\pi_1 X$ to $\Z_2
= \{\,\pm 1\,\}$. As such, it defines a flat complex line bundle $L_{\chi}$.
Since $\chi$ lifts to an integral homology class, the bundle $L_{\chi}$ is
trivial and hence the bundles $P$ and $P\,\otimes\,L_{\chi}$ are isomorphic.
Then $\chi$ acts on $\M^*(X)$ by assigning to a connection $A$ in $P$ the
connection  $A\,\otimes\,\chi$ in $P\,\otimes\,L_{\chi} = P$ induced by $A$ 
and $\chi$. If one views $A$ as a representation $A: \pi_1 X \to SU(2)$, the
above action is given by the formula $\chi(A)(g) = \chi(g) A(g)$ for any $g
\in \pi_1 X$.

\begin{proposition}\label{P:1/2}
The above action of $H^1(X;\Z_2)$ on $\M^*(X)$ is free.
\end{proposition}

\begin{proof}
Let us view $\M^*(X)$ as the irreducible part of the $SU(2)$ representation
variety of $\pi_1 X$. Suppose that $A:\pi_1 X\to SU(2)$ is a fixed point of
$\chi: \M^*(X) \to \M^*(X)$ then there exists a $u \in SU(2)$ such that
$\chi(g) A(g) = u\,A(g)\,u^{-1}$ for all $g\in \pi_1 X$. In particular, by
applying $\chi$ twice, we see that $u^2$ must commute with the image of $A$
in $SU(2)$. Since $A$ is irreducible, this is only possible if $u^2 = \pm 1$.
The case of $u^2 = 1$ should be excluded because then $u = \pm 1$ and $-A(g)
= A(g)$ at least for one $g \in \pi_1 X$, which is impossible. Therefore,
$u^2 = -1$ and, up to conjugation, $u = i$. This means that $A$ is a binary
dihedral representation, that is, its image is contained in $S_i\cup j\cdot
S_i$, where $S_i$ is the circle of unit complex numbers in $SU(2)$.

The representation $A$ maps the subgroup $\pi_1\tilde X = [\pi_1 X,\pi_1 X]$
of $\pi_1 X$ into the commutator subgroup of the binary dihedral group.  The
latter is the unit complex circle; in particular, it is abelian.  Therefore,
the restriction of $A$ onto $\pi_1\tilde X$ factors through $H_1(\tilde X;\Z)
= 0$ and hence $A$ itself factors through $H_1(X;\Z) = \pi_1 X/\pi_1 \tilde X
= \Z$. This contradicts the irreducibility of $A$.
\end{proof}

Since the action of $H^1 (X;\Z_2)$ on $\M^*(X)$ is free,  there exists an
equivariant admissible perturbation $\sigma$ such that $\M^*_{\sigma}(X)$
is non-degenerate and still admits a free action of $H^1(X;\Z_2)$, compare
with Section \ref{S:tori-pert}. This action is orientation preserving, see
Donaldson \cite{donaldson:orientation}, therefore, $\lambda_{FO}(X)$ is at
worst a half--integer. In Section \ref{S:floer}, we will interpret 
$\lambda_{FO}(\overline W)$ as a Floer Lefschetz number and show that 
$\lambda_{FO}(\overline W)$ is always an integer.

%%%%%%%%%%%%%%%%%%%%%%%%%%%%%%%%%%%%%%%%%%%%%%%%%%%%%%%%%%%%%%%%%%%%%%%%%%%

\subsection{Applications}\label{S:fo-appl}
The main topological applications of the Furuta--Ohta invariant stem from
the following two conjectures.

\begin{conjecture}\label{conj1}
Let $X$ be a $\ZZ$--homology $S^1\times S^3$ then $\lambda_{FO}(X) =
\rhofo(X)\pmod 2$.
\end{conjecture}

\begin{conjecture}\label{conj2}
Let $X$ be a $\ZZ$--homology $S^1\times S^3$ which admits an orientation
reversing diffeomorphism inducing an orientation preserving map on
$H^1 (X; \mathbb R)$. Then $\lambda_{FO}(X) = 0$.
\end{conjecture}

First, we show how these two conjectures would disprove the triangulation
conjecture in dimensions $n\ge 5$. Recall from Section \ref{S:rohlin-tri}
that, in order to do that, it would be sufficient to show vanishing of the
Rohlin invariant of any homology sphere $\Sigma$ having second order in
$\Theta^3$. Given such a $\Sigma$, consider a homology cobordism $W$ from
$-\Sigma$ to $\Sigma$ so that $\p W = \Sigma\cup \Sigma$. Identify the
boundary components of $W$ using the identity map. Let $X$ be an orientable
double cover of the resulting non-orientable manifold, then $X$ is a
homology $S^1\times S^3$ admitting an orientation reversing involution
which induces an identity map on $H^1(X;\mathbb R)$. Conjecture \ref{conj2} 
now implies that $\lambda_{FO}(X) = 0$, and Conjecture \ref{conj1} implies 
that $\rho(\Sigma) = \rhofo(X) = \lambda_{FO}(X) = 0\pmod 2$.

That no two homology spheres with non-trivial Rohlin invariant can be 
homology cobordant to each other via a simply connected homology cobordism
$W$ would follow by applying Conjecture \ref{conj1} to the double of $W$. 

%%%%%%%%%%%%%%%%%%%%%%%%%%%%%%%%%%%%%%%%%%%%%%%%%%%%%%%%%%%%%%%%%%%%%%%%%%%%%%

\subsection{Rohlin's invariant and the homotopy $S^1\times S^3$}
\label{S:surgery-4}
Another application of Conjecture \ref{conj1} is to the surgery--theoretic 
classification of smooth manifolds of the homotopy type of $S^1\times S^3$. 
In summary, if surgery theory `worked' in dimension four as it does in 
higher dimensions, then there should exist a fake homotopy $S^1\times S^3$ 
with non-trivial Rohlin invariant. On the other hand, 
Conjecture \ref{conj1} would mean that there
is no such manifold. This would imply a failure of exactness of the smooth
surgery sequence, namely that the $L$--group $L_5(\Z[\Z])$ does not act on
the structure set of $S^1\times S^3$.  Although the existence of
non-diffeomorphic, $s$--cobordant 4--manifolds implies that the surgery
sequence is not exact for simply connected 4--manifolds, we know of no 
example where the group $L_5(\Z[\pi])$ fails to act.

We will briefly review the surgery calculation, and refer the reader
to Wall's book~\cite{wall:book} and the excellent survey of Kirby and
Taylor \cite{kirby-taylor:surgery} for further details and references.
We remind the reader that we are concerned here with {\em smooth}
manifolds; the tools of surgery theory work better in the topological case,
and the calculations are somewhat different. In particular, any homotopy
$S^1 \times S^3$ is homeomorphic to the real one.  The  fake smooth $S^1
\times S^3$ predicted by the surgery sequence would be in a sense the
simplest possible orientable manifold; the fake $\mathbf{R}P^4$ constructed
by Cappell-Shaneson~\cite{cappell-shaneson:rp4} (see
also~\cite{fs:involution}) is also detected by a codimension--one Rohlin
invariant.

To understand this prediction, consider the (hypothetical) surgery
sequence for the structure set of $Y=S^1 \times S^3$,
$$
\begin{CD}
[\Sigma Y,G/PL] @>\theta>> L_5(\Z[\Z]) @>\gamma>> \S(Y) @>N>> [Y,G/PL].
\end{CD}
$$
The hypothetical part on which we concentrate is whether the map called
$\gamma$ giving the action of $L_5(\Z[\Z])$ on $\S(Y) $ is actually
defined. (Kirby and Taylor \cite{kirby-taylor:surgery} explain that a
stabilized (with respect to repeated connected sum with $S^2 \times S^2$)
version of $\gamma$ is defined and fits into an exact surgery sequence
for a `stable' structure set $\bar{\S}$). Since $L_5(\Z) = 0$, there are 
isomorphisms 
\[
L_5(\Z[\Z]) \cong L_5(\Z)\oplus L_4(\Z) 
\stackrel{\sigma/8}{\longrightarrow} \Z
\]
given by a codimension--one signature, see Shaneson \cite{shaneson:product}.  
On the other hand, the calculation
$$
[\Sigma Y,G/PL] = [S^2 \vee S^4 \vee S^5, G/PL] \cong \Z_2 \oplus \Z
$$
and the fact that $\theta$ is a homomorphism, see Wall \cite{wall:book}, 
reduce the calculation of $\theta$ to understanding the {\em
$4$-dimensional} surgery map 
\[
\pi_4(G/PL) \to L_4(\Z)\stackrel{\sigma/8}{\longrightarrow} \Z.
\]    
But an element of $\pi_4(G/PL)$ is a normal map $V^4 \to S^4$, so that the 
normal bundle of $V$ pulls back from a bundle over $S^4$. This implies that
$V$ is in fact spin, so that Rohlin's theorem tells us that the image of 
$\theta$ has index $2$.

In other words, if there were a realization of the action of
$L_5(\Z[\Z])$ on $\S(Y)$ (i.e. if the map $\gamma$ existed) then there
would exist a nontrivial element in $\S(Y)$.  It is easy to check that
there is no self--homotopy equivalence of $Y$ realizing this element.
Hence the conjectured equality of Furuta--Ohta and Rohlin invariants
implies that the surgery sequence is not exact.

The existence of a fake $S^1 \times S^3$ is of course related to the
problem discussed in Section~\ref{S:rohlin-sc} of finding a simply
connected homology cobordism between a Rohlin invariant--one homology
sphere and itself. For gluing the boundary components of such a 
cobordism would give a homotopy $S^1 \times S^3$ with non-trivial Rohlin 
invariant. The converse does not necessarily hold, because a homotopy 
$S^1\times S^3$ would not necessarily have a homology sphere carrying the 
third homology. Similarly, a homology sphere $\Sigma$ with $\rho(\Sigma) 
= 1$ that is of order two in $\Theta^3$ (in the strong sense that $\Sigma 
\# \Sigma$ bounds a contractible manifold) would give rise to a fake
non-orientable manifold homotopy equivalent to $S^1 \widetilde\times
S^3$.  As in the orientable case, the smooth surgery sequence would
predict the existence of such a manifold coming from the action of
$L_5(\Z[\Z^-])$.  See Akbulut \cite{akbulut:fake} for a stabilized
version of this manifold.

%%%%%%%%%%%%%%%%%%%%%%%%%%%%%%%%%%%%%%%%%%%%%%%%%%%%%%%%%%%%%%%%%%%%%%%%%

\section{The Floer homology}\label{S:floer}
Floer \cite{floer:instanton} associated with every homology sphere $\Sigma$
eight abelian groups $I_k(\Sigma)$, $0\le k\le 7$, called (instanton) Floer
homology, which ramify the Casson invariant in that
\[
\lambda(\Sigma) = \frac 1 2\,\sum_k\; (-1)^k \rk I_k(\Sigma).
\]
The definition is as follows. The Floer homology is the homology of the Floer
chain complex $IC_*(\Sigma)$. The free abelian group $IC_k (\Sigma)$ is
generated by the points $A\in \R^*(\Sigma)$ of Floer index $\mu(A) = k\pmod
8$, where the flat moduli space $\R^*(\Sigma)$ may need to be perturbed
first, using admissible perturbations of Section \ref{S:casson-pert}, to
make it non-degenerate. The differential
\[
\p: IC_k (\Sigma) \to IC_{k-1}(\Sigma)
\]
is given by counting (perturbed) ASD connections over the cylinder $\mathbb R
\times \Sigma$ with proper boundary conditions. The Floer homology is
functorial with respect to cobordisms between homology spheres. We will be
mostly interested in the case of a homology cobordism $W$ from a homology
sphere $\Sigma_0$ to a homology sphere $\Sigma_1$. In this case, we have a
well defined homomorphism $W_*: I_* (\Sigma_0) \to I_* (\Sigma_1)$ of degree
zero obtained by counting ASD connections on a trivial $SU(2)$ bundle over
$W$. Again, the ASD equation may need to be perturbed, in which case we use
perturbations of the type described in Section \ref{S:fo-pert}, extended to
match perturbations on $\p W$.

Let $W$ be a homology cobordism from a homology sphere $\Sigma$ to itself
and $W_k: I_k(\Sigma)\to I_k(\Sigma)$ the automorphisms induced by $W$ in
Floer homology. Define the (Floer) Lefschetz number of $W$ by the formula
\[
\Lef(W) = \sum_{k=0}^7\;(-1)^k\;\tr(W_k).
\]
Let $\overline W$ be the closure of $W$ then the usual gluing arguments can
be used to show that
\begin{equation}\label{E:lef=fo}
\Lef(W) = 2\,\lambda_{FO} (\overline W).
\end{equation}

For example, the product cobordism $[0,1] \times \Sigma$ induces the identity
map in Floer homology and hence $\Lef([0,1]\times\Sigma)= 2\lambda(\Sigma)$;
on the other hand, the closure of the product cobordism is $S^1\times\Sigma$,
therefore, $\lambda_{FO} (S^1\times\Sigma) = \lambda (\Sigma)$.  An extension
of this calculation to mapping tori will be discussed in Section \ref{S:map}.

\begin{proposition}
For any homology cobordism $W$ from a homology sphere $\Sigma$ to itself, the
Lefschetz number $\Lef(W)$ is even.
\end{proposition}

\begin{proof}
The homomorphism $W_*$ commutes with the $u$--map of Fr\o yshov
\cite{froyshov}. The result now follows because the $u$--map provides
isomorphisms (over the rationals) $I_k(\Sigma) = I_{k+4}(\Sigma)$ for all $k$.
\end{proof}

Together with Proposition \ref{P:1/2} this implies that the Furuta--Ohta
invariant of the $\ZZ$--homology $S^1\times S^3$ of the type $X = \overline W$
is always an integer. For manifolds of this type, Conjectures \ref{conj1} and
\ref{conj2} can be reformulated as follows. The first conjecture asserts
that $1/2\,\Lef(W) = \rhofo(X) \pmod 2$ or, equivalently, that the modulo
2 reduction of $1/2\,\Lef(W)$ is independent of the choice of $W$. The
second conjecture asserts that $\Lef(W) = 0$ if $W$ admits an orientation
reversing diffeomorphism. 

Since all $\ZZ$--homology $S^1 \times S^3$  that arise in applications 
described in Section \ref{S:fo-appl} are of the type $X = \overline W$, 
proving the above two conjectures about $\Lef(W)$ would suffice for those 
applications. However, it is far from clear that an arbitrary 
$\ZZ$--homology $S^1\times S^3$ (or even a homotopy $S^1\times S^3$ as 
discussed in Section~\ref{S:surgery-4}) can be built in this fashion 
hence the full strength Conjectures \ref{conj1} and \ref{conj2} may be 
needed for applications described in Section \ref{S:surgery-4}. 

%%%%%%%%%%%%%%%%%%%%%%%%%%%%%%%%%%%%%%%%%%%%%%%%%%%%%%%%%%%%%%%%%%%%%%%%%%%%%%

\section{Mapping tori}\label{S:map}
Let $\Sigma$ be a homology sphere and $\tau: \Sigma\to\Sigma$ an orientation
preserving diffeomorphism. Recall that the mapping torus of $\tau$ is the
smooth 4-manifold $X_{\tau} = ([0,1]\times\Sigma)\,/(0,x)\sim (1,\tau(x))$
with the product orientation. In this section, we give explicit formulas for
$\lambda_{FO}(X_{\tau})$ in the case of finite order $\tau$; note that the
theory has a rather different character depending on whether $\tau$ has
fixed points or not. Our formulas first express $\lambda_{FO}(X_{\tau})$ in
terms of the equivariant Casson invariant, and then identify the latter as a
linear combination of the (regular) Casson invariant and certain classical
knot invariants. We use these explicit formulas to verify Conjectures
\ref{conj1} and \ref{conj2} for the mapping tori of finite order
diffeomorphisms.

%%%%%%%%%%%%%%%%%%%%%%%%%%%%%%%%%%%%%%%%%%%%%%%%%%%%%%%%%%%%%%%%%%%%%%%%%%%%%%

\subsection{The definition of the equivariant Casson invariant}
Let $\tau: \Sigma \to \Sigma$ be an orientation preserving diffeomorphism of
finite order. It induces a map $\tau^*: \R^*(\Sigma) \to \R^*(\Sigma)$ via
pull back of flat connections. Let $\R^{\tau}(\Sigma)$ be the fixed point
set of $\tau^*$.  After perturbing the flatness equation $F_A = 0$, if 
necessary, $\R^{\tau}(\Sigma)$ is a compact canonically oriented manifold of 
dimension zero. We define the equivariant Casson invariant as
\[
\lambda^{\tau}(\Sigma) = \frac 1 2\,\#\R^{\tau}(\Sigma),
\]
where $\#\R^{\tau}(\Sigma)$ stands for the signed count of points in 
$\R^{\tau}(\Sigma)$.

As we see, the definition follows closely that of the (regular) Casson
invariant; however, there are a few important differences. First, the
non-degeneracy of $\R^{\tau}(\Sigma)$ at a point $A$ means vanishing of the
equivariant cohomology group $H^1_{\tau}(\Sigma;\ad A)$, which may be
strictly smaller than the group $H^1(\Sigma;\ad A)$. Therefore, we can
use equivariant perturbations as in \cite{CS} to achieve non-degeneracy
of $\R^{\tau}(\Sigma)$ without achieving non-degeneracy of the full moduli
space $\R^*(\Sigma)$ (the latter should not be expected anyway because of
the equivariant transversality problem). The equivariant perturbations in
question are first constructed on the quotient manifold $\Sigma/\tau$ by
applying the procedure of Section \ref{S:casson-pert} to a link $\gamma_k$
in $\Sigma/\tau$, and then lifted to $\Sigma$.

Second, the space $\R^{\tau}(\Sigma)$ is oriented using the equivariant
spectral flow, which is the spectral flow of the operators (\ref{E:ka})
restricted to the subspace of $(\Omega^0\,\oplus\,\Omega^1)(\Sigma,\ad P)$
invariant with respect to the induced action of $\tau$. This spectral flow
may well be different from the spectral flow used to orient $\R^*(\Sigma)$;
in short, the natural inclusion $\R^{\tau}(\Sigma)\to \R^*(\Sigma)$ is not
necessarily orientation preserving.

%%%%%%%%%%%%%%%%%%%%%%%%%%%%%%%%%%%%%%%%%%%%%%%%%%%%%%%%%%%%%%%%%%%%%%%%%%%%%%

\subsection{Furuta--Ohta vs. equivariant Casson}
This section is dedicated to expressing the Furuta--Ohta invariant for mapping
tori of finite order diffeomorphisms in terms of the equivariant Casson
invariant. We give a brief outline of the proof here and refer the reader
to our paper \cite{ruberman-saveliev:mappingtori} for all the details.

\begin{theorem}\label{T:fo}
If $\tau: \Sigma \to \Sigma$ has finite order then $\lambda_{FO}(X_{\tau})
= \lambda^{\tau}(\Sigma)$.
\end{theorem}

\begin{proof}
First observe that all solutions of the equation $F^+_A = 0$ on $X_{\tau}$
are flat by Chern--Weil theory. A flat connection over $X_{\tau}$ is pulled
back to a flat connection over $\Sigma$ via an inclusion $i: \Sigma \to
X_{\tau}$, and this pull back map defines a two-to-one map $\M^*(X_{\tau})
\to \R^{\tau}(\Sigma)$. The latter can be easily seen by interpreting flat
connections as representations of respective fundamental groups. We have
a splitting exact sequence
\[
\begin{CD}
0 @>>> \pi_1(\Sigma) @>>> \pi_1(X_{\tau}) @>>> \Z @>>> 0.
\end{CD}
\]
Let $t$ be a generator of $\Z$ then every irreducible representation $A:
\pi_1 (X_{\tau}) \to SU(2)$ determines and is uniquely determined by the
pair $(\alpha,u)$ where $u = A(t)$ and $\alpha = i^* A: \pi_1 \Sigma \to
SU(2)$ is an irreducible representation such that $\tau^*\alpha = u\alpha
u^{-1}$ (in particular, the conjugacy class of $\alpha$ belongs to
$\R^{\tau}(\Sigma)$). Replacing $u$ by $-u$ in the above gives rise to a
new representation $\pi_1 (X_{\tau}) \to SU(2)$, which in fact is simply
the image of $A$ under the action of $H^1 (X;\Z_2)$ on $\M^*(X_{\tau})$.
This results in a two-to-one correspondence between $\M^*(X_{\tau})$ and
$\R^{\tau}(\Sigma)$. Moreover, a direct calculation with cohomology shows
that $H^1(X_{\tau};\ad A) = H^1_{\tau}(\Sigma; \ad i^* A)$.

Let us assume for the moment that $\R^{\tau}(\Sigma)$ is non-degenerate.
Then so is $\M^*(X_{\tau})$, and to prove the identity $\lambda_{FO}
(X_{\tau}) = \lambda^{\tau}(\Sigma)$, we only need to show that the above
identification of the moduli spaces is orientation preserving. This can
be achieved by first interpreting the orientation on $\M^*(X_{\tau})$ in
terms of the orientation transport of a family of ASD operators $D_{A(t)}$,
see Nicolaescu \cite{nicolaescu:seiberg-witten}. Next, the orientation
transport can
be identified with the orientation given by the equivariant spectral flow
along the path $i^*A(t)$ by expanding forms on $X_{\tau}$ into Fourier
series in the direction of $S^1$ and comparing the spectra of operators
$D_{A(t)}$ and $K_{i^*A(t)}$, compare with Atiyah, Patodi and Singer
\cite{aps:I}.

Another way to compare the orientations is by viewing $X_{\tau}$ as an
(orbifold) $S^1$--bundle over the quotient manifold $\Sigma/\tau$. One
can apply adiabatic limit techniques to this situation; after that the
result will follow from the fact that the equivariant spectral flow on
$\Sigma$ equals the (regular) spectral flow on the quotient $\Sigma/\tau$.

Finally, if $\R^{\tau}(\Sigma)$ fails to be non-degenerate, perturb the
flatness equation $F_A = 0$ into $F_A = *\nabla h$ where $h$ is lifted
from a perturbation $h'$ on $\Sigma/\tau$. According to \cite{CS}, there
are enough such perturbations to make $\R^{\tau}(\Sigma)$ non-degenerate.
Next, perturb the ASD equation $F^+_A = 0$ into $F^+_A = \sigma(A)$ where
$\sigma(A)$ is the self--dual part of the 2--form $*\nabla h'$ pulled back
to $X_{\tau}$ via the projection $X_{\tau}\to\Sigma/\tau$. Note that
$\sigma(A)$ is of the type described in Section \ref{S:fo-pert} with 
$N^3_k = S^1\times D^2$. One can verify that there still exists the 
two-to-one correspondence $\M^*_{\sigma}(X_{\tau})\to\R^{\tau}_h(\Sigma)$ 
between the perturbed moduli spaces, and that the perturbations $\sigma$ 
as described above are sufficient to make $\M^*_{\sigma}(X_{\tau})$ 
non-degenerate. The latter essentially follows by comparing the abundancy 
concepts for the two types of perturbations and using the equality 
$H^1(X_{\tau};\ad A) = H^1_{\tau}(\Sigma;\ad i^* A)$ of the respective 
Zariski tangent spaces.
\end{proof}

%%%%%%%%%%%%%%%%%%%%%%%%%%%%%%%%%%%%%%%%%%%%%%%%%%%%%%%%%%%%%%%%%%%%%%%%%%

\subsection{Equivariant Casson: non-free actions}\label{S:non-free}
Let $\tau: \Sigma \to \Sigma$ be an orientation preserving diffeomorphism
of finite order $n$,  and suppose that the fixed point set of $\tau$ is
non-empty. Then the quotient manifold $\Sigma'=\Sigma/\tau$ is a homology
sphere, and the projection $\Sigma \to \Sigma'$ is a branched covering
with branch set a knot $k \subset \Sigma'$. The following is proved in 
Collin--Saveliev \cite{CS}.

\begin{theorem}\label{T:non-free}
In the situation described above, the equivariant Casson invariant is given
by the formula
\[
\lambda^{\tau}(\Sigma) = n\cdot\lambda(\Sigma') + \frac 1 8\,
\sum_{m=0}^{n-1}\; \sign^{m/n} (k),
\]
where $\sign^a (k)$ is the Tristram--Levine equivariant knot signature,
defined as the signature of the Hermitian form $(1 - e^{2\pi i a}) S +
(1 - e^{-2\pi i a}) S^t$, for any choice of Seifert matrix $S$ of $k$.
\end{theorem}

\begin{proof}
The proof proceeds by pushing equivariant flat connections from $\Sigma$
down to singular connections over the quotient $\Sigma'$ and using Herald's
theorem \cite{herald:signatures} which asserts that an appropriate count
of the latter connections is a certain linear combination of 
$\lambda(\Sigma')$ and equivariant knot signatures.
\end{proof}

\begin{corollary}\label{C:non-free}
Let $\tau: \Sigma \to \Sigma$ be a finite order orientation preserving
diffeomorphism having fixed points then $\lambda_{FO}(X_{\tau}) =
\rhofo(X_{\tau})\pmod 2$.
\end{corollary}

\begin{proof}
This is equivalent to showing that $\lambda^{\tau}(\Sigma) = \rho(\Sigma)
\pmod 2$ which follows from Theorem \ref{T:non-free} by standard 
techniques of geometric topology, see for instance Viro \cite{viro}.
\end{proof}

%%%%%%%%%%%%%%%%%%%%%%%%%%%%%%%%%%%%%%%%%%%%%%%%%%%%%%%%%%%%%%%%%%%%%%%%%%

\subsection{Equivariant Casson: free actions}\label{S:free}
In case $\tau$ acts freely on $\Sigma$,  the quotient $\Sigma' =
\Sigma/\tau$ is a homology lens space. It is easy to see that $\Sigma'$
can be obtained by $(n/q)$--surgery on a knot $k$ in a homology sphere
$Y$ where $n$ is the order of $\tau$ and $q$ is relatively prime to $n$.
The following is proved in our paper \cite{ruberman-saveliev:mappingtori}.

\begin{theorem}\label{T:free}
In the situation described above, the equivariant Casson invariant is given
by the formula
\[
\lambda^{\tau}(\Sigma) = n\cdot\lambda(Y) + \frac 1 8\,\sum_{m=0}^{n-1}\;
\sign^{m/n} (k) + \frac q 2\,\Delta''_k(1),
\]
where $\Delta_k(t)$ is the Alexander polynomial of $k\subset Y$ normalized
so that $\Delta_k(1) = 1$ and $\Delta_k(t) = \Delta_k(t^{-1})$.
\end{theorem}

\begin{proof}
The proof proceeds by pushing equivariant flat connections from $\Sigma$
to the homology lens space $\Sigma'$ and identifying the count of the
latter connections with a sum of invariants of the type discussed by 
Boyer--Nicas \cite{boyer-nicas} and Boyer--Lines \cite{boyer-lines}. Note 
that these invariants are different from Walker's invariant 
\cite{walker:casson} in that they only count \emph{irreducible} flat 
connections over $\Sigma'$; their gauge theoretic definition is implicit in 
work of Cappell, Lee and Miller \cite{clm:III}. An application of surgery 
formulas and Herald's theorem on equivariant knot signatures 
\cite{herald:signatures} completes the proof.
\end{proof}

Let $Y_n$ be the $n$--fold cyclic cover of $Y$ branched along $k$. Observe
that $Y_n$ is a homology sphere and that $\Sigma$ is obtained by
$(1/q)$--surgery on $Y_n$ along a lift $k_n$ of $k$. Combining Theorem
\ref{T:free} with Theorem \ref{T:non-free} we see that
\[
\lambda^{\tau}(\Sigma) = \lambda^{\tau}(Y_n) + \frac q 2\,\Delta''_k(1).
\]
Since we already know that $\lambda^{\tau}(Y_n) = \rho(Y_n)\pmod 2$,
see Corollary \ref{C:non-free}, we will be able to conclude that
$\lambda^{\tau}(\Sigma) = \rho(\Sigma)\pmod 2$ once we verify the
following result, see \cite{ruberman-saveliev:mappingtori}.

\begin{lemma}
Let $Y$ be an integral homology sphere and $\pi: Y_n\to Y$ its $n$--fold
cyclic branched covering with branch set a knot $k$. Let $k_n$ be the knot
$\pi^{-1}(k)$ in $Y_n$. If $Y_n$ is an integral homology sphere then
$\arf(k_n) = \arf(k)\pmod 2$.
\end{lemma}

\begin{corollary}\label{C:free}
Let $\tau: \Sigma \to \Sigma$ be a free finite order orientation preserving
diffeomorphism then $\lambda_{FO}(X_{\tau}) = \rhofo(X_{\tau})\pmod 2$.
\end{corollary}

%%%%%%%%%%%%%%%%%%%%%%%%%%%%%%%%%%%%%%%%%%%%%%%%%%%%%%%%%%%%%%%%%%%%%%%%%%%

\subsection{Orientation reversal}
Corollaries \ref{C:non-free} and \ref{C:free} verify Conjecture \ref{conj1}
for the mapping tori $X_{\tau}$ of finite order diffeomorphisms $\tau:
\Sigma \to \Sigma$. To verify Conjecture \ref{conj2} for such mapping tori,
observe that any orientation reversing diffeomorphism $f: X_{\tau} \to
X_{\tau}$ which preserves homology orientation can be viewed as a
diffeomorphism $f: -X_{\tau}\to X_{\tau}$ preserving both orientation and
homology orientation, so that $\lambda_{FO}(-X_{\tau}) = \lambda_{FO}
(X_{\tau})$. On the other hand, since $-X_{\tau}$ is the mapping torus of
the orientation preserving diffeomorphism $\tau: -\Sigma \to -\Sigma$, we
conclude that $\lambda_{FO}(-X_{\tau}) = \lambda^{\tau}(-\Sigma)$ by Theorem
\ref{T:fo}. The equivariant Casson invariant changes sign with the change of
orientation, therefore, $\lambda_{FO}(-X_{\tau}) = -\lambda^{\tau}(\Sigma) =
-\lambda_{FO}(X_{\tau})$ and thus $\lambda_{FO}(X_{\tau})$ must vanish.

%%%%%%%%%%%%%%%%%%%%%%%%%%%%%%%%%%%%%%%%%%%%%%%%%%%%%%%%%%%%%%%%%%%%%%%%%%%

\section{Examples}\label{S:ex}
In this section, we give examples of finite order maps $\tau: \Sigma \to
\Sigma$ for which the Furuta--Ohta invariant and the Floer Lefschetz number
can be computed explicitly using methods described above. An interesting
observation is that the map $(W_{\tau})_*: I_*(\Sigma) \to I_*(\Sigma)$
induced by the mapping cylinder $W_{\tau}$ of $\tau$ need not be identity
even when $\tau$ acts trivially on the representation variety $\R^*(\Sigma)$.
This has to do with the fact that $(W_{\tau})_*$ is defined using a count of
ASD connections over $W_{\tau}$ with signs which are not determined solely by
the fundamental group. 

%%%%%%%%%%%%%%%%%%%%%%%%%%%%%%%%%%%%%%%%%%%%%%%%%%%%%%%%%%%%%%%%%%%%%%%%%%%

\subsection{Akbulut cork}\label{S:ex-cork}
\begin{figure}[!ht]
\begin{center}
\psfrag{0}{$0$}
\includegraphics{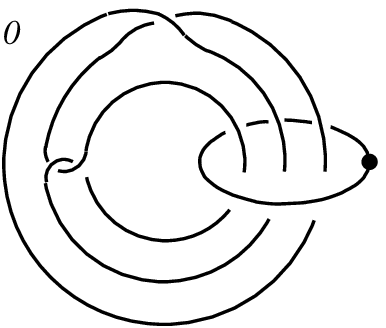}
\end{center}
\caption{}
\label{cork}
\end{figure}

By Akbulut cork we mean the smooth contractible 4--manifold $W$ obtained by
attaching a two-handle to $S^1\times D^3$ along its boundary as shown in
Figure \ref{cork}. It can be embedded into a blown up elliptic surface $E(n)
\#(-\mathbb C P^2)$ in such a way that cutting it out and re-gluing by an
involution on $\Sigma = \p W$ changes the smooth structure on $E(n)\#
(-\mathbb C P^2)$ but preserves its homeomorphism type, see Akbulut
\cite{akbulut:cork} and Gompf and Stipsicz \cite{gompf}.

\begin{figure}[!ht]
\begin{center}
         \begin{minipage}[b]{0.38\linewidth}
           \begin{center}
           \psfrag{0}{$0$}
           \includegraphics{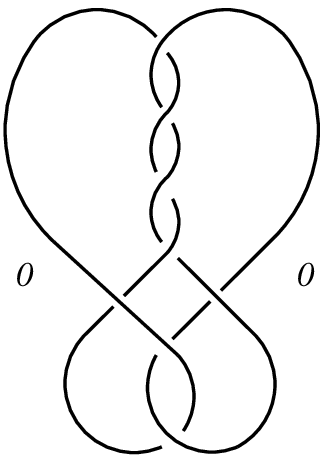}
           \end{center}
           \caption{}\label{symmcork}
          \end{minipage}
          \begin{minipage}[b]{0.38\linewidth}
           \begin{center}
           \psfrag{k}{$k^*$}
           \includegraphics{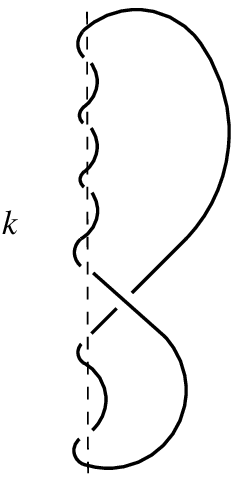}
           \end{center}
           \caption{}\label{corkquotient}
          \end{minipage}
\end{center}
\end{figure}

The involution $\tau:\Sigma \to \Sigma$ simply interchanges the two link
components, $k_1$ and $k_2$; this is best seen when the link is drawn in
a symmetric form as in Figure \ref{symmcork}.

The manifold $\Sigma' = \Sigma/\tau$ is obtained from $S^3$ by surgery on
the knot $k^*$ which is the image of the link $k_1 \cup k_2$,  see Figure
\ref{corkquotient}. Note that the canonical longitudes of $k_1$ and $k_2$
project onto a longitude of $k^*$ whose linking number with $k^*$ equals
one. The dotted line in Figure \ref{corkquotient} represents the branch
set $k$ of $\Sigma \to \Sigma'$. This picture can be viewed as a surgery
description of the knot $k$ in $\Sigma' = S^3$. A little exercise in Kirby
calculus shows that $k$ can be obtained from the left--handed
$(5,6)$--torus knot on six strings by adding one full left--handed twist
on two adjacent strings. The signature of $k$ can differ by at most two
from the signature of the left--handed $(5,6)$--torus knot, which equals
16. Since $\sign k$ must be divisible by eight, we conclude that $\sign k
= 16$. Therefore, $\lambda_{FO}(X_{\tau}) = \lambda^{\tau}(\Sigma) = 2$
by Theorem \ref{T:non-free} and $\Lef(W_{\tau})= 2\cdot\lambda_{FO}(X_{\tau})
= 4$, see (\ref{E:lef=fo}).

In fact, the latter formula provides enough information to describe the map
$(W_{\tau})_*: I_*(\Sigma)\to I_*(\Sigma)$. According to \cite{saveliev:cork},
the Floer homology of $\Sigma$ is trivial in even degrees, and is a copy of
$\Z$ in each of the odd degrees. Therefore, $\tau_*: I_k (\Sigma) \to I_k
(\Sigma)$ is necessarily plus or minus identity for each $k$. Since the
Lefschetz number of $\tau_*$ equals $4$, this implies that $(W_{\tau})_*$ is
minus identity.

It is worth mentioning that the character variety of $\Sigma$ is
non-degenerate, and that the action induced by $\tau$ on $\R^*(\Sigma)$ is
a non-trivial permutation, see \cite{saveliev:cork}.

%%%%%%%%%%%%%%%%%%%%%%%%%%%%%%%%%%%%%%%%%%%%%%%%%%%%%%%%%%%%%%%%%%%%%%%%%%

\subsection{Seifert fibered homology spheres}
Let $\Sigma(a_1,\ldots,a_n)$  be a Seifert fibe\-red homology sphere viewed
as a link of singularity, and $\tau$ the involution on it
induced by complex conjugation. It makes $\Sigma(a_1,\ldots,a_n)$ into a
double branched cover of $S^3$ with branch set a Montesinos knot $k$. A
straightforward application of Theorem \ref{T:non-free}
shows that both $\lambda_{FO}(X_{\tau})$ and $(1/2)\,\Lef(W_{\tau})$ are
equal to $1/8\,\sign k$, which is also known as the $\bar\mu$--invariant of
$\Sigma(a_1,\ldots,a_n)$, see Neumann \cite{neumann} and Siebenmann
\cite{siebenmann}. As a side note, we mention that the homology cobordism
invariance of the $\bar\mu$--invariant, shown with the help of the orbifold
Seiberg--Witten theory, is the main tool in proving Theorem \ref{T:ffs}.

As a rule, the $\bar\mu$--invariant of $\Sigma(a_1,\ldots,a_n)$ differs 
from its Casson invariant,  which indicates that the action $\tau_*: 
I_*(\Sigma(a_1,\ldots,a_n))\to I_*(\Sigma(a_1,\ldots,a_n))$ is 
non-trivial. With an extra effort, this action can be described 
explicitly. According to Fintushel--Stern \cite{fs:instanton}, the group 
$I_k (\Sigma(a_1,\ldots,a_n))$ is trivial if $k$ is even, and is a free 
abelian group of rank, say, $b_k$ if $k$ is odd. The following is proved 
in \cite{saveliev:brieskorn}.

\begin{theorem}
The map $\tau_*: I_k (\Sigma(a_1,\ldots,a_n)) \to I_k (\Sigma(a_1,
\ldots,a_n))$ is identity if $k = 1\pmod 4$ and minus identity if
$k = -1\pmod 4$.
\end{theorem}

The character variety of any Seifert fibered homology sphere with (no more
than) three singular fibers is non-degenerate, and one can show that the
action induced on it by $\tau$ is trivial, see \cite{saveliev:brieskorn}.
The same is true for any number of singular fibers after one perturbs the
character variety using equivariant perturbations as in
\cite{saveliev:seifert}.

%%%%%%%%%%%%%%%%%%%%%%%%%%%%%%%%%%%%%%%%%%%%%%%%%%%%%%%%%%%%%%%%%%%%%%%%%%

\section{Circle bundles}\label{S:circle}
In this section we describe another class of $\ZZ$--homology $S^1\times S^3$
for which both Conjectures \ref{conj1} and \ref{conj2} have been verified;
it turns out that in all of these examples, both invariants $\lambda_{FO}$
and $\rhofo$ vanish.  As mentioned in Section \ref{S:fo-ex}, these
examples were pointed out to us by Scott Baldridge.

%%%%%%%%%%%%%%%%%%%%%%%%%%%%%%%%%%%%%%%%%%%%%%%%%%%%%%%%%%%%%%%%%%%%%%%%%%

\subsection{The manifolds}
Let $Y$ be a closed oriented 3--manifold having integral homology of $S^1
\times S^2$, and denote by $\Delta(t)$ its Alexander polynomial normalized
so that $\Delta(1) = 1$ and $\Delta(t^{-1}) = \Delta(t)$. Every manifold
$Y$ as above can be obtained by 0--surgery on a knot $k$ in a homology
sphere $\Sigma$; the Alexander polynomial of $Y$ then equals the Alexander
polynomial of $k\subset\Sigma$.

Let $\pi: X\to Y$ be an $S^1$--bundle classified by its Euler class $e\in
H^2(Y;\Z) = \Z$ and assume that $e = 1$.  Then one can use the Gysin
exact sequence to show that $X$ has integral homology of $S^1\times S^3$.
To ensure that $X$ is a $\ZZ$--homology $S^1\times S^3$, the universal
abelian cover $\tilde X$ should have homology of $S^3$; this is equivalent
to having $\Delta(t) = 1$. These are the manifolds $X$ that we will study 
in this section.

%%%%%%%%%%%%%%%%%%%%%%%%%%%%%%%%%%%%%%%%%%%%%%%%%%%%%%%%%%%%%%%%%%%%%%%%%%%%

\subsection{Calculating $\rhofo(X)$}
Recall that  $\rhofo (X) = \rho (M, \sigma)$  where $M \subset X$ is an
embedded 3--manifold representing a generator in $H_3(X;\Z)$ and $\sigma$
is the induced spin structure. A calculation with Gysin exact sequence
shows that $M$ can be obtained by taking a surface $S$ carrying $H_2(Y;\Z)$
and then taking its preimage in $X$, which is of course a circle bundle over
$S$. The spin structure on $X$ can be described as follows. There is an exact
sequence of bundles
\[
\begin{CD}
0 @>>> V @>>> TX @>>> \pi^*\,TY @>>> 0
\end{CD}
\]
where $V$ is the vertical tangent space. This induces a bijection between
spin structures on $Y$ and those on $X$. Thus the induced spin structure
$\sigma$ on $M$ is given by taking a spin structure on $Y$ (either one will
do) and restricting to $S$, and then using the same exact sequence to get
a spin structure on $M$.

To calculate the Rohlin invariant $\rho(M,\sigma)$, consider the disk bundle
$W \to S$ with Euler class 1, whose boundary is $M$. The manifold $W$ is not
spin but we still obtain
\[
\rho(M,\sigma) = (\sign W - S\cdot S)/8 + \arf(S) = \arf(S)\pmod 2,
\]
where $\arf(S)$ is the Arf--invariant of the induced spin structure. Since
the Alexander polynomial of $Y$ vanishes, we conclude that $\arf(S) = 0$.

%%%%%%%%%%%%%%%%%%%%%%%%%%%%%%%%%%%%%%%%%%%%%%%%%%%%%%%%%%%%%%%%%%%%%%%%%%

\subsection{Calculating $\lambda_{FO}(X)$}
In this section, we find it convenient to work with $SO(3)$ flat connections
rather than with $SU(2)$ ones. The two settings are equivalent due to the
presence of a free orientation preserving $H^1(X;\Z_2)$ action on $\M^*(X)$
whose quotient is the $SO(3)$ flat moduli space $\M^*(X,SO(3))$. Note that
this remains true after a small generic equivariant perturbation $\sigma$
making $\M^*_{\sigma}(X)$ non-degenerate, see Proposition \ref{P:1/2} and
Section \ref{S:tori-pert}.

Without loss of generality we will assume that $\pi_2(Y) = 0$. The homotopy
exact sequence of the $S^1$--bundle $\pi: X \to Y$ then implies that $\pi_1
X$ is a central extension of $\pi_1 Y$ by the integers,
\[
\begin{CD}
1 @>>> \Z @>>> \pi_1 X @>\pi_* >> \pi_1 Y @>>> 1.
\end{CD}
\]
Let $h$ be a generator in $\Z$ then every irreducible representation $A:
\pi_1 X \to SO(3)$ has the property that $A(h) = 1$. Therefore, we have
a natural identification $\pi^*: \R^*(Y,SO(3)) \to \M^*(X,SO(3))$ where
$\R^*(Y,SO(3))$ is the $SO(3)$ character variety of $\pi_1 Y$.

This identification induces an isomorphism $\pi^*: H^1(Y;\ad A) \to
H^1(X;\ad\pi^* A)$ of Zariski tangent spaces, which is easily seen from
the Gysin exact sequence
\[
\begin{CD}
0 @>>> H^1(Y;\ad A) @>\pi^* >> H^1(X;\ad\pi^* A) @>>> H^0(Y;\ad A) @>>> 0
\end{CD}
\]
after one notes that $H^0(Y;\ad A)$ vanishes because $A$ is irreducible.
Thus the moduli space $\M^*(X,SO(3))$ is non-degenerate if and only if
$\R^*(Y,SO(3))$ is. Should the non-degeneracy fail, both $\R^*(Y,SO(3))$
and $\M^*(X,SO(3))$ need to be perturbed. This can be done in a consistent
manner as in the proof of Theorem \ref{T:fo} so that the pull back of 
connections via $\pi: X \to Y$ still provides a bijective correspondence
between the perturbed moduli spaces. To show that this correspondence is
orientation preserving, one can use adiabatic limit techniques, see
Nicolaescu \cite{nicolaescu:private}.

Thus $\lambda_{FO}(X)$  equals one half times the signed count of points
in the (perturbed) moduli space $\R^*(Y,SO(3))$. In computing the latter,
we will rely on the papers \cite{herald:signatures} and
\cite{ruberman-saveliev:mappingtori}, to which we refer the reader for all
the details.

View $Y$ as the result of 0--surgery on a knot $k$ in a homology sphere
$\Sigma$. Let $Z = \Sigma\setminus N(k)$ be the knot $k$ exterior, and
$m$ and $\ell$ be the canonical meridian and longitude on the torus
$\p Z$. Then the dual torus $\mathcal P = \R(\p Z,U(1))$ has coordinates
$(\phi,\psi)$ such that the holonomies along $m$ and $\ell$ are equal to
$\exp(i\phi)$ and $\exp(i\psi)$, respectively. The inclusion $\p Z\to Z$
induces a natural restriction map from a double cover of $\R(Z,SU(2))$ to
$\mathcal P$ whose image, $\tilde\C$, is generically an immersed curve.
Moreover, since $\Delta_k (t) = 1$, this image consists of finitely many
circles.

Let us consider the splitting $\R^*(Y,SO(3)) = \R^*_0(Y,SO(3))\,\cup
\,\R^*_1(Y,SO(3))$, where $\R^*_w (Y,SO(3))$ stands for the moduli space
of irreducible flat connections on an $SO(3)$ bundle $P$ having  $w_2(P)
= w\in H^2(Y;\Z_2) = \Z_2$. Up to an overall constant, the signed count
of points in $\R^*_1(Y,SO(3))$ equals the intersection number of
$\tilde\C$ with the circle $\psi = \pi$ and hence equals $\Delta''_k(1)
= 0$. Similarly, the signed count of points in $\R^*_0 (Y,SO(3))$ equals,
up to an overall constant, the intersection number of $\tilde\C$ with
the circle $\psi = 0$. Since the circles $\psi = 0$ and $\psi = \pi$ are
homologous in $\mathcal P$, this intersection number is again equal to
$\Delta''_k(1) = 0$.

%%%%%%%%%%%%%%%%%%%%%%%%%%%%%%%%%%%%%%%%%%%%%%%%%%%%%%%%%%%%%%%%%%%%%%%%%%%%

\section{A general approach to the conjectures}\label{S:surgery}
To prove Conjecture \ref{conj1}, we would like to carry out a surgery
program similar to the program showing the equality $\lambda(\Sigma) =
\rho(\Sigma)\pmod 2$ in three dimensions. Many essential features of that
proof, like reliance on Heegaard splittings and the Alexander polynomial,
are not easily transferable one dimension higher. Therefore, our first
step is giving a purely gauge theoretic proof of $\lambda(\Sigma) =
\rho(\Sigma)\pmod 2$, see \cite{ruberman-saveliev:casson}. This proof is
briefly described in Section \ref{S:surgery-3}. After that, we outline a
surgery program in dimension four for proving Conjecture \ref{conj1}. The
final step of this program, the calculation of the degree zero Donaldson
polynomial for homology tori, is now complete (modulo some technical
issues), see Section \ref{S:tori},
while intermediate steps require a lot of work.

Verifying Conjecture \ref{conj2} amounts to comparing the Furuta--Ohta
invariants for a $\ZZ$--homology $S^1\times S^3$ with two opposite
orientations. We only know how to prove this conjecture in the situation
when the ASD moduli space $\M^*(X)$ is non-degenerate; the general case
is still out of our reach.

%%%%%%%%%%%%%%%%%%%%%%%%%%%%%%%%%%%%%%%%%%%%%%%%%%%%%%%%%%%%%%%%%%%%%%%%%%%

\subsection{Homology 3--tori and the Casson invariant}\label{S:surgery-3}
The Casson invariant has the property that $\lambda(\Sigma) = \rho(\Sigma)
\pmod 2$ for any homology sphere $\Sigma$, see Section \ref{S:casson-prop}.
The traditional proof of this property involves several steps.  

The first step takes us back to the original definition of the Casson 
invariant in terms of Heegaard splittings. This step relies on Taubes' 
theorem \cite{taubes:casson}.  

In the second step, surgery formulas are used to express the Rohlin and 
the Casson invariants in terms of their first, second and third 
difference quotients, see \cite{akbulut-mccarthy} or \cite{saveliev:casson}. 
The third difference quotients, $\rho'''(\Sigma)$ and $\lambda'''(Y)$, are
invariants of a manifold $Y$ having integral homology of the 3--torus (or a 
homology 3--torus, for short). Since the surgery formula for the Rohlin 
invariant happens to be the modulo two reduction of the surgery formula for 
the Casson invariant, showing that $\lambda(\Sigma) = \rho(\Sigma)\pmod 2$ 
amounts to verifying that $\lambda'''(Y) = \rho'''(Y)\pmod 2$. 

According to Kaplan \cite{kaplan:even}, the invariant $\rho'''(Y)$ equals 
the determinant of $Y$ defined as $\det Y = (a_1\cup a_2\cup a_3)\,[Y]
\pmod 2$, where $a_1$, $a_2$, $a_3$ is a basis in $H^1(Y;\Z_2)$. The final 
step is then proving that $\lambda'''(Y) = \det Y\pmod 2$. This is 
accomplished by using the interpretation of Casson's difference quotients 
in terms of the Alexander polynomials of links.

Our approach eliminates the first step of the above program and hence the 
need for Taubes' theorem. Instead, we directly apply the Casson surgery 
formula, which in gauge theoretic terms is an easy corollary of the Floer 
exact triangle, see \cite{floer:knots} and \cite{braam-donaldson:knots}. 
This expresses $\lambda(\Sigma)$ in terms of third difference quotients 
$\lambda'''(Y)$ where $\lambda'''(Y)$ now equals a proper count of flat
connections on an $SO(3)$ bundle over $Y$ with a non-trivial second
Stiefel--Whitney class $w$ (the invariant $\lambda'''(Y)$ does not depend
on the choice of $w\ne 0$).

The final step of the proof, that is, showing that $\lambda'''(Y) = \det Y
\pmod 2$, is done gauge--theoretically along the lines of our proof of
Theorem \ref{T:tori}. We refer the reader to \cite{ruberman-saveliev:casson}
for all the details.

%%%%%%%%%%%%%%%%%%%%%%%%%%%%%%%%%%%%%%%%%%%%%%%%%%%%%%%%%%%%%%%%%%%%%%%%%%%%%

\subsection{Round surgery}\label{S:surgery-round}
We would like to carry out a similar program in dimension four. The idea is
to follow the above proof step by step in a category of 4--manifolds $X$
over $S^1$, and reduce Conjecture \ref{conj1} to the result about
$\ZZ$--homology 4--tori proved in Theorem \ref{T:tori}.

More precisely, a 4--manifold over $S^1$ is a smooth closed oriented
manifold $X$ of dimension four with a preferred cohomology class $\alpha
\in H^1(X;\Z)$. There are obvious notions of cobordism and surgeries of
manifolds over $S^1$.    Roughly speaking, we want to think of our
4--manifolds as `looking like' they are the product of $S^1$ with a
3--manifold, and thus the surgeries allowed should look like they
are really given by $S^1$ times a surgery on the 3--manifold. This brings
into our discussion the concept of round surgery defined by Asimov
\cite{asimov:round}. A four--dimensional round handle of index $k$ is a
pair
\[
S^1 \times (D^k \times D^{3-k}, S^{k-1} \times D^{3-k})
\]
attached to a manifold with boundary. There are obvious notions of round
handle decompositions, round surgeries etc. Asimov \cite{asimov:round}
showed that the existence (in dimensions at least four) of round handle
decompositions is governed by the Euler characteristic. From this it is
not difficult to prove the following result.

\begin{proposition}\label{P:round}
Let $X$ be a smooth $\ZZ$--homology $S^1\times S^3$. Then $X$ may be
obtained from $S^1 \times S^3$ via a series of round surgeries of indices
$2$ and $3$.
\end{proposition}

Unfortunately, the intermediate stages in the resulting cobordism might
not have well defined Rohlin invariants (and their gauge theory invariants
are not so good either).   The problem would be rectified if we could make
the intermediate stages look like $S^1$ times a 3--manifold. To accomplish
this, we need to restrict the kind of surgeries we allow.  The restriction
is in terms of the {\em degree} of the round surgery, which is by
definition the integer $\langle \alpha, S^1 \rangle$.   If the degree of a
round surgery is $\pm 1$, then we can define appropriate gauge theoretic
and Rohlin invariants, and then try to compare them.  Thus we need to
strengthen Proposition~\ref{P:round} to say that $X$ can be obtained
by a series of round surgeries, each of which has degree $\pm 1$.

This strengthening, although it sounds fairly innocent, seems to be very
difficult to accomplish, and indeed there may need to be some modification
in the program which we now describe. Assuming that there is some
strengthened proposition along these lines,  here are some ideas on how we
might obtain a proof of Conjecture~\ref{conj1}.   Except for the very last
one, each of the steps in this outline still has many details to be filled
in.

The main point is that there seem to exist surgery formulas for both
Rohlin and Furuta--Ohta invariants for round surgeries of degree $\pm 1$.
These have a form that one might expect by thinking about surgery on $S^1$
times a 3--manifold.  Namely, different framings for a surgery on $X$ give
rise to manifolds $X_1$ and $X_0$, where $X_1$ has the same homology as
$X$, and $b_1(X_0) = b_1(X) +1$.  (The notations are supposed to suggest
$+1$ and $0$ surgery on a 3--manifold, respectively). Both the Rohlin and
the Furuta--Ohta invariants of $X$ can be written as the sum of invariants
for $X_0$ and $X_1$.   As in the 3--dimensional theory, we need to have a
direct verification that the two invariants coincide for manifolds with
sufficiently large $b_1$; it turns out that $b_1 = 4$ is as far as we need
to go.  We have actually carried out this final step;  in the next section
we define the Rohlin and Furuta--Ohta invariants for homology 4--tori and
outline our theorem that they in fact coincide.

%%%%%%%%%%%%%%%%%%%%%%%%%%%%%%%%%%%%%%%%%%%%%%%%%%%%%%%%%%%%%%%%%%%%%%%%%%

\section{Gauge theory on homology tori}\label{S:tori}
We study the degree zero Donaldson polynomial for $\ZZ$--homology 4--tori
and relate it to a properly defined Rohlin invariant. A similar result
holds for homology 3--tori, which leads to a purely gauge theoretic proof 
of the fact  that the Casson invariant of a homology sphere reduces modulo 
2 to its Rohlin invariant, compare with Section \ref{S:surgery-3}.

%%%%%%%%%%%%%%%%%%%%%%%%%%%%%%%%%%%%%%%%%%%%%%%%%%%%%%%%%%%%%%%%%%%%%%%%%%

\subsection{$\ZZ$--homology 4--tori}\label{S:tori-def}
By $\ZZ$--homology 4--torus we mean a closed oriented smooth spin
4--manifold $X$ such that $H_*(X;\Z) = H_*(T^4;\Z)$ and $H_*(\tilde X_a;\Z)
\allowbreak = H_*(T^3;\Z)$, where $\tilde X_a$ is the infinite cyclic cover
of $X$ corresponding to a choice of primitive element $a \in H^1(X;\Z)$.
The intersection form of $X$ on the second cohomology is always isomorphic
to the sum of three copies of the hyperbolic 2--form; however,
the cup--product on the first cohomology of $X$ may vary.  Let $a_0$, $a_1$,
$a_2$, and $a_3$ be a basis in $H^1(X;\Z_2)$ then $\det X = (a_0\cup a_1
\cup a_2\cup a_3)\,[X]\pmod 2$ is independent of the choice of $a_0$, $a_1$,
$a_2$, and $a_3$ and is called the determinant of $X$. A $\ZZ$--homology
4--torus $X$ is called odd if $\det X = 1\pmod 2$, and is called even
otherwise.

\begin{figure}[!ht]
\begin{center}
\psfrag{0}{$0$}
\psfrag{SS}{$3\,\# (S^1\times S^2)$}
\psfrag{T}{$T^3$}
\includegraphics[totalheight=1.1in]{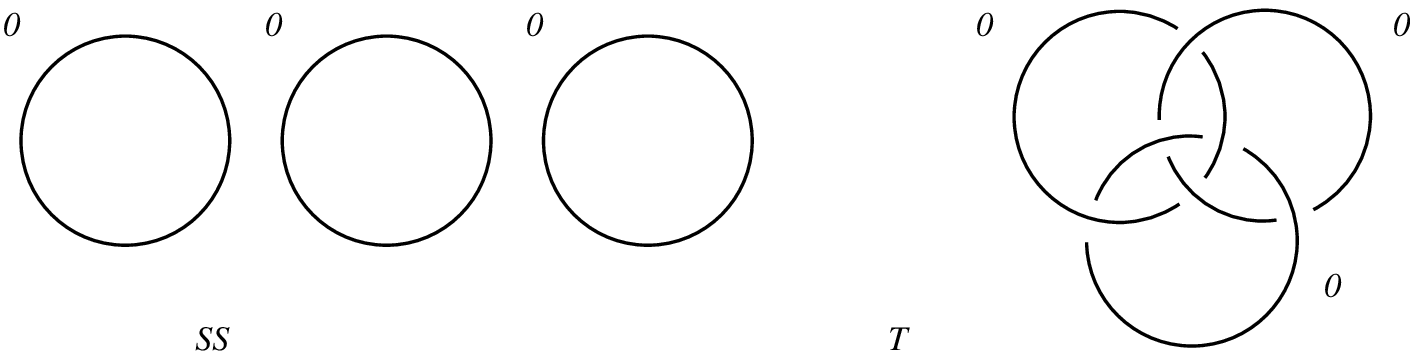}
\end{center}
\caption{}
\label{F:tori}
\end{figure}

There are any number of ways to construct $\ZZ$--homology 4--tori;
for instance, any homology 3--torus times a
circle is a $\ZZ$--homology 4--torus.  In turn, an ample source of homology
3--tori is surgery on any 3--component link with zero framing matrix as in
Figure \ref{F:tori}. The $\ZZ$--homology 4--torus $S^1 \times (3\,\# (S^1
\times S^2))$ is even while the torus $T^4 = S^1 \times T^3$ is odd.

\bigskip

\begin{figure}[!ht]
\begin{center}
\psfrag{P}{$D_a$}
\psfrag{Q}{$D_b$}
\psfrag{0}{$0$}
\psfrag{p}{$-q$}
\psfrag{q}{$q$}
\includegraphics[totalheight=1.5in]{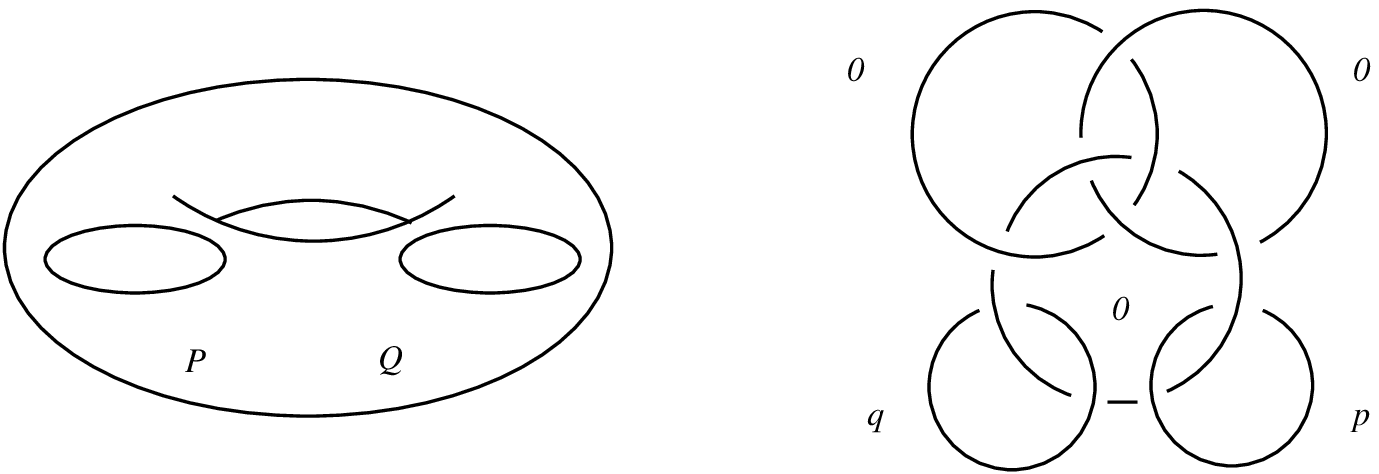}
\end{center}
\caption{}
\label{F:log}
\end{figure}

Another interesting family of examples of $\ZZ$--homology 4--tori is obtained
by the following construction.  View $T^4$ as a trivial $T^2$--bundle over
$T^2$ and consider two embedded disks $D_a$ and $D_b$ in the base $T^2$, see
Figure \ref{F:log}. 

\noindent
Cut out $D_a \times T^2$ and $D_b\times T^2$ and glue them back using 
automorphisms
\begin{alignat*}{2}
(a,x,y) &\to (a^q\, x^{q-1},a x, y) & &\quad\text{on $\p D_a\times T^2$},  \\
(b,x,y) &\to (b^{-q}\, x^{-q-1},b x, y)& & \quad \text{on $\p D_b\times T^2$}.
\end{alignat*}
If $q$ is odd, the resulting manifold is an odd $\ZZ$--homology 4--torus which
we denote by $T^4(q,-q)$. It is a double branched cover of a log--transform of
the Kummer surface with exceptional curves blown down. In fact, $T^4 (q,-q)
= S^1 \times T^3(q,-q)$, where $T^3(q,-q)$ is the surgery on the link shown
in Figure \ref{F:log}.

It is not necessarily the case that a $\ZZ$--homology 4--torus is a
product with $S^1$.  Examples of this can arise from the round
surgery construction discussed in Section \ref{S:surgery-round}. Gluing
round 2--handles along a collection of embedded disjoint 2--tori, if done
right, will result in a $\ZZ$--homology 4--torus.

%%%%%%%%%%%%%%%%%%%%%%%%%%%%%%%%%%%%%%%%%%%%%%%%%%%%%%%%%%%%%%%%%%%%%%%%%%%

\subsection{The Rohlin invariant}\label{S:tori-rohlin}
Let $X$ be a $\ZZ$--homology 4--torus with spin structure $\sigma$, and let
$a \in H^1(X;\Z)$ be a primitive element such that the cyclic cover
$\tilde X_a \to X$ corresponding to $a: \pi_1 X \to \Z$ has integral homology
of $T^3$. Let $M \subset X$ be Poincar\'e dual to $a$ and define the Rohlin
invariant
\[
\rho(X,a,\sigma) = \rho(M,\sigma)\pmod 2
\]
using the induced spin structure $\sigma$ on $M$. Complete $a$ to a basis 
$a, x_1, x_2, x_3 \in H^1(X;\Z_2)$ and define
\[
\bar\rho(X,a) = \sum_{x\in\Span\{x_1, x_2, x_3\}}\rho(X,a,\sigma + x).
\]

\begin{theorem} The invariant $\bar\rho(X,a)$ is well defined and $\bar\rho
(X,a) = \det X\pmod 2$.
\end{theorem}

We refer the reader to our paper \cite{ruberman-saveliev:donaldson} for a
complete proof of this theorem.

%%%%%%%%%%%%%%%%%%%%%%%%%%%%%%%%%%%%%%%%%%%%%%%%%%%%%%%%%%%%%%%%%%%%%%%%%%%

\subsection{The Donaldson invariant}\label{S:tori-don}
Let $X$ be a $\ZZ$--homology 4--torus, $P \to X$ an $SO(3)$--bundle with
$p_1(P) = 0$ and $w_2 (P)\ne 0$, and $\G_{SU(2)}$ the group of automorphisms
of $P$ that lift to $SU(2)$. The space $\M(P)$ of ASD connections on $P$
modulo $\G_{SU(2)}$ has formal dimension $-2p_1(P)\,[X] - 3\,(1 - b_1 + b^2_+)
(X) = 0$. The condition $w_2(P) \ne 0$ implies that all connections in $\M(P)$
are irreducible, and it follows from the Chern--Weil theory that all of them
are actually flat.

After perturbing the ASD equation $F^+_A = 0$ using a generic admissible
perturbation $\sigma$ as described in Section \ref{S:fo-pert} we obtain a
perturbed moduli space $\M_{\sigma}(P)$ which is a compact oriented manifold
of dimension zero, compare with Donaldson \cite{donaldson:orientation}.  The
degree zero Donaldson polynomial is defined as
\[
\D_0 (X,P) = \#\,\M_{\sigma}(P).
\]
An alternative approach to the definition of $\D_0(X,P)$, which we employ in
\cite{ruberman-saveliev:donaldson}, is via $U(2)$ bundles and projectively
ASD connections with fixed central part.

In order to state our main result about $\D_0 (X,P)$,  we need to impose two
extra hypotheses. One is a mild technical restriction on the bundle $P$\,: 
we assume that there exists $\xi \in H^1(X;\Z_2)$ such that $w_2(P) \cup \xi 
\ne 0$ (note that this condition is automatically satisfied for odd 
$\ZZ$--homology 4--tori). The other hypothesis is needed because, at present, 
we have some unresolved issues in completing the proof of the existence
of equivariant perturbations, see Section \ref{S:tori-pert}. We assume that
\begin{itemize}
     \item[(*)]  \emph{there is an equivariant generic admissible
perturbation that makes the ASD moduli space into a smooth
$0$-manifold.}
    \end{itemize}
We are quite confident that the ideas briefly sketched in
Section~\ref{S:tori-pert} will show that hypothesis (*) is always satisfied.

\begin{theorem}\label{T:tori}
Under the above restriction on the bundle $P$ (and, for the moment, hypothesis 
(*)), the invariant $\D_0 (X,P)$ is divisible by four and
\[
\frac 1 4 \; \D_0 (X,P) = \det X \pmod 2.
\]
\end{theorem}

Note that the formula of Theorem \ref{T:tori} need not hold over the
integers. For instance, $\D_0 (T^4(q,-q),P)$ is equal to $\pm\, 4 q^2$ if
$P$ is a pull back from $T^3(q,-q)$, and to $\pm\, 4$ otherwise, see
\cite{ruberman-saveliev:donaldson}. It is also worth mentioning that, 
in the spirit of Witten's conjecture relating Donaldson and
Seiberg--Witten invariants, Theorem \ref{T:tori} is consistent with the
result of the first author and S.~Strle \cite{ruberman-strle:tori} on the
(mod 2) evaluation of the Seiberg--Witten invariant for homology
4--tori.  The main result of this paper is that the Seiberg--Witten
invariant of the $\text{Spin}^c$ structure associated to a spin
structure on a homology 4--torus $X$ is congruent to $\det X$ modulo 2.

%%%%%%%%%%%%%%%%%%%%%%%%%%%%%%%%%%%%%%%%%%%%%%%%%%%%%%%%%%%%%%%%%%%%%%%%%%%%

\subsection{A sketch of the proof of Theorem \ref{T:tori}}\label{S:tori-sk}
Let us assume for the moment that the moduli space $\M(P)$ is non-degenerate.
Let $\G_{SO(3)}$ be the full gauge group of automorphisms of $P$, then we
have the exact sequence
\[
1 \longrightarrow \G_{SU(2)} \longrightarrow \G_{SO(3)} \longrightarrow
H^1 (X;\Z_2) \longrightarrow 1
\]
so that the moduli space $\M(P)$ is acted upon by $H^1 (X; \Z_2)$ and its
quotient is the moduli space $\overline\M(P)$ of $\G_{SO(3)}$ equivalence
classes of ASD connections on $P$. This action is a complete analogue of
the action described in Section \ref{S:fo-prop}.

The action of  $H^1(X;\Z_2) = (\Z_2)^4$  is orientation preserving,  see
Donaldson \cite{donaldson:orientation}, therefore, in order to compute
$\D_0(X,P)$, one can count not individual points in $\M(P)$ but rather
their orbits. From this point on, we will proceed by showing that there
are no orbits of orders one or two, so that $(1/4)\cdot \D_0(X,P)\pmod 2$
equals the number of the 4--orbits (because the eight-- and the
sixteen--orbits do not contribute to the above count).

If the moduli space $\M(P)$ fails to be non-degenerate, it is perturbed
using an admissible perturbation $\sigma$ which is {\it equivariant} with
respect to the action of $H^1(X;\Z_2)$.  Achieving non-degeneracy of 
$\M_{\sigma}(P)$ using equivariant perturbations is not an easy task, see 
Section \ref{S:tori-pert}, but once it is done, the above counting 
argument can be applied to $\M_{\sigma}(P)$ to complete the proof.

%%%%%%%%%%%%%%%%%%%%%%%%%%%%%%%%%%%%%%%%%%%%%%%%%%%%%%%%%%%%%%%%%%%%%%%%%%%%%

\subsection{The holonomy correspondence}\label{S:tori-hol}
To do an actual count of the four--orbits, we take advantage of the fact
that the connections in $\M(P)$ are flat and hence can be interpreted
algebraically using the holonomy map. More precisely, let $w_2(P) = w \in
H^2(X;\Z_2)$ then the usual holonomy correspondence identifies the moduli
space $\overline{\M}(P)$ with a compact subset of the character variety 
$\R(X,SO(3))$, which we call $\R_w(X,SO(3))$. The proper algebraic tool 
for lifting this correspondence to the covering $\M(P)\to\overline\M(P)$ 
is projective representations.

A map  $\rho: \pi_1 X \to SU(2)$  is called a projective representation if
$\rho(gh) = c(g,h)\allowbreak\rho(g)\rho(h)$ where $c(g,h)$ belongs to
$\Z_2=\{\pm 1\}$ viewed as the center of $SU(2)$. The 2--cocycle $c$
defines a cohomology class $[c] \in H^2 (\pi_1 X; \Z_2)$ such that $[c] =
w_2(\ad\rho)$. Here, $w_2 (\ad\rho)$ stands for the second Stiefel--Whitney
class of the flat $SO(3)$ bundle with holonomy $\ad\rho$, and the above
equality makes sense after we identify $H^2(\pi_1 X;\Z_2)$ as a natural
subset of $H^2(X;\Z_2)$. In general, $H^2(\pi_1 X;\Z_2)$ is not equal to
$H^2 (X;\Z_2)$\,;\; this minor point needs a separate treatment, see
\cite{ruberman-saveliev:donaldson}.

Define $\P_c\,(X;SU(2))$ to be the set of conjugacy classes of projective
representations $\rho$ with fixed $c$. It is acted upon by $H^1 (X;\Z_2)$ 
by the rule $(a,\rho)\to a\cdot\rho$. According to 
\cite{ruberman-saveliev:casson}, the holonomy correspondence 
$\overline\M(P)\to \R_w(X,SO(3))$ lifts to a bijective correspondence 
$\M(P) \to \P_c\,(X, SU(2))$ so that we have the following commutative 
diagram
\[
\begin{CD}
\M(P) @>= >> \P_c\,(X, SU(2)) \\
@VVV  @V\ad VV \\
\overline{\M}(P) @>= >> \R_w\,(X,SO(3))
\end{CD}
\]
where both vertical arrows are the quotient maps by the action of
$H^1(X;\Z_2)$. This interpretation of $\M(P)$ in terms of projective
representations allows for a precise description of the four--orbits.

%%%%%%%%%%%%%%%%%%%%%%%%%%%%%%%%%%%%%%%%%%%%%%%%%%%%%%%%%%%%%%%%%%%%%%%%%%%%

\subsection{The four--orbits}\label{S:tori-four}
Let $a, b\in H^1(X;\Z_2)$ be two different non-trivial elements stabilizing
$\rho \in \P_c\,(X, SU(2))$. The argument from the proof of Proposition
\ref{P:1/2} when applied simultaneously to $a$ and $b$ shows that, possibly
after conjugation, $\im\rho\subset (S_i\cup j\cdot S_i)\cap (S_j\cup k\cdot
S_j) = \{\pm 1,\pm i,\pm j,\pm k\}$.  Equivalently, $\im\ad\rho\subset \Z_2
\oplus \Z_2 \subset SO(3)$ where the latter inclusion is given by $(A,B)\to
A\oplus B\oplus A\cdot B$. The same argument can further be used to show
that the action of $H^1(X;\Z_2)$ on $\P_c\,(X, SU(2))$ does not have orbits
of orders one or two.

Thus the set of the four--orbits is in a bijective correspondence with the set
of representations $\ad\rho:\pi_1 X\to SO(3)$ which factor through a subgroup
$\Z_2\,\oplus\,\Z_2$ of $SO(3)$, modulo $SO(3)$ conjugation. Every such a
representation induces representations $\alpha, \beta: \pi_1 X \to \Z_2$,
which can be viewed as cohomology classes $\alpha, \beta \in H^1 (X;\Z_2)$.
This identifies the set of the 4--orbits with  $\Lambda^2_0\;H^1(X;\Z_2)$,
the set of decomposable elements in the second exterior power of
$H^1(X;\Z_2)$. A straightforward calculation gives $w_2(\ad\rho) = \alpha
\cup\beta \in H^2 (X;\Z_2)$, hence the set of the four--orbits with fixed
$w$ can be identified with the preimage of $w$ under the map
\[
\cup: \Lambda^2_0\; H^1 (X;\Z_2) \longrightarrow H^2 (X;\Z_2).
\]
It is an exercise in algebraic topology to show that, under the hypotheses
of Theorem \ref{T:tori}, the above map is a bijection if $\det X$ is odd
and is zero if $\det X$ is even.

%%%%%%%%%%%%%%%%%%%%%%%%%%%%%%%%%%%%%%%%%%%%%%%%%%%%%%%%%%%%%%%%%%%%%%%%%%%%%

\subsection{Perturbations}\label{S:tori-pert}
The above argument only applies as stated if the moduli space $\M (P)$ is
non-degenerate, otherwise, $\M(P)$ needs to be perturbed first. It turns 
out that the four--orbits are always non-degenerate and remain such after 
a small enough perturbation, so the argument of Section \ref{S:tori-four}
counting the four--orbits goes unchanged. Thus the eight-- and
sixteen--orbits are the only ones that are perturbed.   As of this
writing, some details of what is briefly described below remain to be
completely verified.

We define our perturbations $\sigma$ using holonomy around loops $\psi_k$
as in Section \ref{S:fo-pert}.  In the equivariant setting, we have to
address two new issues.  First, we need to make $\sigma$ equivariant with
respect to the $H^1 (X;\Z_2)$ action; this is done by requiring that $0 =
[\psi_k] \in H_1 (X;\Z_2)$. Next, we need to show that these equivariant
perturbations are generic, that is, there are enough of them to make
$\M_{\sigma}(P)$ non-degenerate. Perturbations along disjoint loops as
in Section \ref{S:fo-pert} are sufficient if the action of $H^1(X;\Z_2)$
is free; otherwise, one may need more general perturbations. We refer the 
reader to the forthcoming revision of \cite{ruberman-saveliev:donaldson}
for all the details.

It is worth mentioning that, since the quotient space of $\M(P)$ by the
$H^1(X;\Z_2)$ action is the $SO(3)$ moduli space $\overline\M (P)$, the
equivariant perturbation theory described above is essentially the
$SO(3)$ perturbation theory.

%%%%%%%%%%%%%%%%%%%%%%%%%%%%%%%%%%%%%%%%%%%%%%%%%%%%%%%%%%%%%%%%%%%%%%%%%%%%%

\end{document}